\documentclass[12pt]{article}

\def\real{\mathbb{R}}

\def\Sd{\mbox{{\boldmath{${\mathbb S}$}}$^d$}}

\newcommand{\sd}{\mathbb{S}}

\newcommand{\T}{\mathcal{T}}

\newcommand{\Div}{\mbox{\rm Div\,}}

\newcommand{\skalar}[1]{\left\langle #1 \right\rangle}

\newcommand{\dual}[2]{\left\langle #1 \right\rangle_{{#2}^*\times #2}}

\newtheorem{Theorem}{\sc Theorem}

\newtheorem{Lemma}[Theorem]{\sc Lemma}
\newtheorem{Corollary}[Theorem]{\sc Corollary}

\def\sqr#1#2{{
    \vcenter{
         \vbox{\hrule height.#2pt
               \hbox{\vrule width.#2pt height#1pt \kern#1pt
                     \vrule width.#2pt
               }
               \hrule height.#2pt
         }
    }
}}

\def\div{\mathop{\rm div}\nolimits}
\def\Div{\mathop{\rm Div}\nolimits}

\def\bar{\overline}
\def\real{\mathbb{R}}

\newcommand{\R}{{\if mm {\rm I}\mkern -3mu{\rm R}\else \leavevmode
\hbox{I}\kern -.17em\hbox{R} \fi}}

\def\lista#1
{{ \itemindent 0.0cm \labelsep .2cm \leftmargin 0.8cm \rightmargin
0.0cm \labelwidth 0.6cm \topsep 0.0mm
\parsep 0.0mm
\itemsep 0.0mm
\begin{list}{}
{ \setlength{\leftmargin}{1.0cm} \setlength{\rightmargin}{0.0cm}
\setlength{\parsep}{.0mm} \setlength{\topsep}{.0mm}
\setlength{\parskip}{.0cm} \setlength{\itemsep}{.0cm} }
{#1}\end{list}} }

\usepackage{amsmath}
\usepackage{amssymb}
\usepackage[mathscr]{eucal}
\usepackage{graphicx}
\usepackage{epstopdf} 
\usepackage{color}

\usepackage{textcomp}
\usepackage{subfig}
\usepackage{placeins}
\usepackage{url}

\usepackage{latexsym}
\usepackage{mathrsfs}
\usepackage{amsfonts,amsmath,amssymb}
\usepackage{indentfirst}
\usepackage{subeqnarray}
\usepackage{verbatim}

\usepackage{booktabs}
\usepackage{geometry}
\usepackage{bm}
\usepackage{enumerate}
\usepackage{float}
\usepackage{caption}
\begin{document}

\title{ Numerical analysis of a non-clamped dynamic thermoviscoelastic contact problem \ }

\author{
Piotr Bartman$^{\,1}$, Krzysztof Bartosz$^{\,1}$,\\ Micha{\l} Jureczka$^{\,1}$, Pawe{\l} Szafraniec$^{\,1}$,
  \
\\ ~ \\
{\small $^1$ Jagiellonian University, Faculty of Mathematics and Computer Science} \\
{\small ul. Lojasiewicza 6, 30348 Krakow, Poland}
}

\date{}

\maketitle

\thispagestyle{empty}

\begin{center}
	{\it The paper is dedicated to Professor Weimin Han} \\
	{\it on the occasion of his 60th birthday}
\end{center}
\vskip 4mm

\noindent {\bf Abstract.} In this work, we analyze a non-clamped dynamic viscoelastic contact problem involving thermal effect. The friction law is described by a non-monotone relation between the tangential stress and the tangential velocity. This leads to a system of second-order inclusion for displacement and a parabolic equation for temperature. We provide a fully discrete approximation of the problem and find optimal error estimates without any smallness assumption on the data. The theoretical result is illustrated by numerical simulations.   

\vskip 4mm


\vskip 4mm


\newpage

\section{Introduction}\label{Introduction}

Physical contact processes appear in many fields of industry and everyday life. Hence, it is reasonable to study mathematical models concerning various properties and parameters of physical bodies in contact with each other. In particular, thermal effects play a special role in such models. A physical body may change its properties and shape when heated. On the other hand, contact between two physical bodies is directly related to a heat exchange. Additionally, frictional contact usually leads to heat generation. It is clear that a well-formulated mathematical model should take into account all of the above phenomena.

There are numerous publications concerning mathematical analysis of contact processes with thermal effect, which essentially deal with the existence and uniqueness of a solution for the stated problems. We recommend \cite{ADLY, AKRS, A2002, A1997, CHAU2, FT, RS2000} as representative references in this field.

Apart from the theoretical aspect of mathematical modeling, numerical analysis also plays a significant role. Typically, it includes approximation schemes based on temporal and spatial discretization, an estimate of the error between the approximate and exact solutions, and numerical simulations. This kind of approach carried out for various mechanical contact problems with or without thermal effect is presented, for instance, in \cite{ABF, BFV, BFT, CFHS, CFKSE, FM, HSS, HILD_RVNARD, khlare, kpr}.

In this paper we deal with the dynamic contact of a viscoelastic body with a foundation involving a thermal effect. The contact is frictional, and the friction law is modeled by a non-monotone relation between the tangential velocity and the tangential stress. More precisely, the relation has the form of an inclusion involving Clarke subdifferential of a locally Lipschitz non-convex potential. It is worth mentioning that such inclusions are closely related to the concept of hemivariational inequalities (HVIs), which are used as a variational formulation of various non-monotone contact problems of mechanics. The theory and applications of HVIs have been broadly developed in the last few decades. 
For example, we refer to \cite{MOS} for a comprehensive review of recent achievements in this field. The numerical approach to stationary HVIs based on the Finite Element Method (FEM) was first investigated in \cite{Haslinger}. In recent years, similar techniques supported by numerical simulations have been applied to the numerical analysis of various contact processes in mechanics modeled by Clarke subdifferential inclusions or HVIs. We refer to \cite{BBD, Bartosz, BBK, BDS, FHH, HW, WWH} for further details.

The non-monotone character of the friction law introduces one of the main difficulties in our model. Despite this issue, we provide a numerical analysis of the problem based on the fully discrete scheme. The main result consists of the estimate of the error between the discrete and the exact solution. We show that the error estimate is linear with respect to the temporal and spatial discretization parameters in the case where the spatial discretization is based on the first-order FEM. We also present results from numerical simulations that illustrate and validate the linear rate of convergence with respect to the discretization parameters. 

This paper provides a generalization of the result obtained in \cite{BDS}. In contrast to the model studied there, we deal with the non-clamped body, which introduces additional difficulty to the problem. Furthermore, our results are obtained without smallness assumptions, which are required in \cite{BDS}. These advances are possible due to the arguments contained in Corollary \ref{Erhilng_Gamma}.

The remainder of the paper is organized as follows. In Section \ref{Preliminaries} we introduce some preliminary material. In Section \ref{Mechanical_problem_and_variational_formulation} we formulate the mathematical model of the dynamic thermoviscoelastic contact problem and present its variational formulation. The main results on the error estimates for the fully discrete numerical scheme are presented in Section \ref{SecError}. Finally, in Section \ref{SecNum}, we present numerical simulations of a two- and three-dimensional contact problem and provide numerical evidence of optimal order convergence for the linear finite elements.


\section{Mathematical preliminaries}\label{Preliminaries}
In this section, we briefly present the definitions and notation used in the paper. In what follows, we denote by $\langle \cdot, \cdot \rangle_{ X}$ 
the duality pairing between a Banach space 
$X$ and its dual $X^*$. Throughout the paper, we denote by $C$ 
a generic constant whose value may change from line to line.
The symbol $C(0, T; X)$ represents the space of continuous functions on $[0, T]$ with values in $X$. For a function $u\colon[0,T]\to X$
we denote by $\dot u$ and $\ddot u$ its first and second time derivatives, respectively.  

We recall the definitions of the generalized directional 
derivative and the generalized Clarke subgradient for a locally Lipschitz function $\varphi \colon X \to \real$, where $X$ is a Banach space (see~\cite{CLARKE}). 
The generalized directional derivative of $\varphi$ at the point $x \in X$ in direction $v \in X$, denoted by $\varphi^{0}(x; v)$, is defined by
$$\displaystyle 
\varphi^{0}(x; v) =
\limsup_{y \to x, \ t\downarrow 0}
\frac{\varphi(y+tv) - \varphi(y)}{t}.
$$
The generalized subgradient of $\varphi$ at $x$, denoted by
$\partial \varphi(x)$, is a subset of a dual space $X^*$ given by
$$
\partial \varphi(x) = \{ \, 
\zeta \in X^* \mid \varphi^{0}(x; v) \ge 
{\langle \zeta, v \rangle}_{X^* \times X}  
\ \ \mbox{\rm for all} \ v \in X \, \}. 
$$

Let $\real^d$ denote an Euclidean space of dimension $d$, with the scalar product and the norm given by
\begin{align}
u\cdot v& =  u_iv_i,\quad\,\,\,\,
\|v\|_{\real^d} =(v\cdot v)^{1/2}\,\,\ \ \mbox{\rm for all} \
u=(u_i),\,v=(v_i)\in \mathbb{R}^d.\nonumber
\end{align}
Hereafter, all indices $i,j,k,l$ run between $1$ and $d$. Moreover, we apply the so-called summation convention over repeated indices, e.g. the notation $a_ib_i$ means $\sum_{i=1}^da_{i}b_{i}$ and $a_{ij}b_{ij}$ means $\sum_{i,j=1}^da_{ij}b_{ij}$. We denote by $\Sd$ the space of second-order symmetric tensors in $\real^d$, or equivalently the space of symmetric
matrices $d\times d$. The scalar product and the norm in $\Sd$ are given by
\begin{align}
\sigma:\tau& =\sigma_{ij}\tau_{ij},\quad
\|\tau\|_{\sd^d}=(\tau:\tau)^{1/2}\,\,\ \ \mbox{\rm for all} \  \sigma=(\sigma_{ij}),\,
\tau=(\tau_{ij})
\in\mathbb{S}^{d}.\nonumber
\end{align}

\indent Let $\Omega\subset\real^d$ be a bounded domain with a Lipschitz boundary $\Gamma=\partial\Omega$. 
For any function $u\colon\Omega\to\real^d$ such that $u=(u_i)$, $u_i\colon\Omega\to\real$, we use the symbol $u_{i,j}$ for its partial derivative $\frac{\partial u_i}{\partial x_j}$. Furthermore, we use symbols $\div u$ and $\varepsilon(u)$  to denote its divergence operator and the deformation operator, respectively, defined by
\[\div u=u_{i,i},\quad  \varepsilon(u) = (\varepsilon_{ij}(u)),\quad
\varepsilon_{ij}(u) = \frac{1}{2}( {u}_{i,j} + {u}_{j,i} )
. \] 

For a function $\sigma\colon\Omega\to\sd^d$, such that $\sigma=(\sigma_{ij})$, $\sigma_{ij}\colon\Omega\to\real$, we use the symbol $\Div\sigma$ to denote its divergence operator given by
\[ \Div\sigma = (\sigma_{ij,j} ).\]
In the sequel, we will use the following spaces of functions defined in $\Omega$.\\[2mm] 
	$$H=L^2(\Omega;\real^d),\qquad Q=\{\,\sigma=(\sigma_{ij})\,\,|\,\,
	\sigma_{ij}=\sigma_{ji} \in L^{2}(\Omega)\,\}.$$
The spaces $H$ and $Q$ are real Hilbert spaces with inner products given by
	$$( u,v)_H  = \int_{\Omega}{{u}_{i}{v}_{i}\,dx},\qquad 
	( \sigma,\tau )_Q = \int_{\Omega}{\sigma_{ij}\tau_{ij}\,dx. }$$
The norms induced by the scalar products mentioned above are denoted by
$\|\cdot\|_{H}$ and $\|\cdot\|_Q$, respectively.

Let ${\nu}=(\nu_i)$ be the outward unit normal vector at the boundary
$\partial\Omega$.  We denote by
$v_\nu$ and $v_\tau$  the normal and tangential parts of the vector field 
$v\colon\partial\Omega\to\real^d$,
defined by $v_\nu=v\cdot\nu$ and $v_\tau=v-v_\nu\nu$, respectively. Similarly, 
$\sigma_\nu$ and $\sigma_\tau$ represent the normal and tangential parts of the tensor field $\sigma\colon\partial\Omega\to\Sd$ defined by
$\sigma_\nu=(\sigma\nu)\cdot\nu$ and $\sigma_\tau=\sigma-\sigma_\nu\nu$, respectively.

\medskip

We now recall two Green formulas (cf. \cite{MOS}, (2.6) and (2.7), respectively) needed to obtain a variational formulation of the problem of contact mechanics. 
\begin{equation}\label{Green_a}
\int_\Omega(u \div v+\nabla u\cdot v)dx=\int_\Gamma u(v\cdot \nu)\,d\,\Gamma \ \ \mbox{\rm for all} \  u\in W^{1,p}(\Omega), v\in W^{1,q}(\Omega;\real^d),
\end{equation}
where $1\leq p<\infty$, $\frac{1}{p}+\frac{1}{q}=1$ and

\begin{equation}\label{Green}
(\sigma,\varepsilon(v))_Q+ ( \Div \sigma,\,v)_H=\int_\Gamma\sigma\nu\cdot v\,d\,\Gamma \ \ \mbox{\rm for all} \  v\in H^1(\Omega;\real^d),\,\,\sigma\in C^1(\Omega;\Sd).
\end{equation}

\medskip

We complete this section with the following two lemmas, which will be used throughout the paper.

\begin{Lemma}[{\cite{DMP2}, Lemma 8.4.12}]\label{Erhling}
	Let $X$, $Y$ and $Z$ be Banach spaces such that 
	$X$ is compactly embedded in $Y$, 
	and $Y$ is continuously embedded in $Z$. 
	Then, for every $\varepsilon >0$, there exists a constant $C(\varepsilon) > 0$ such that
	\begin{equation}\label{erh}
	\nonumber
	\|x\|_Y\leq \varepsilon \, \|x\|_X 
	+ C(\varepsilon) \, \|x\|_Z 
	\ \ \mbox{\rm for all} \ \ x \in X.
	\end{equation}
\end{Lemma}

\begin{Lemma}[\cite{HS}, Lemma 7.25]\label{gronwall} Let $T>0$ be fixed, $N>1$ and $k=T/N$. Assume $\{e_n\}_{n=1}^N$ and $\{g_n\}_{n=1}^N$ be two sequences of non-negative numbers satisfying
	\begin{equation*} 
	e_n \leq c g_n + c\sum_{j=1}^n ke_j \quad \mbox{for} \ n=1,\ldots,N,
	\end{equation*}
	 with a positive constant $c$ independent of $N$. Then there exists a positive constant $\hat{c}$, independent of $N$, such that
	\begin{equation}
	\max_{1\le n\le N} e_n \le \hat{c} \max_{1\le n \le N} g_n.\nonumber
	\end{equation}
\end{Lemma}

\section{Mechanical problem and variational formulation}\label{Mechanical_problem_and_variational_formulation}
In this section, we describe the mathematical model of the thermoviscoelastic dynamic contact problem in its classical and variational formulations. 

We consider a body that, in its reference configuration, occupies a bounded domain $\Omega \subset \mathbb{R}^d$, $d=2,3$ with a Lipschitz boundary $\Gamma$, consisting of two parts $\Gamma_N$ and $\Gamma_C$, such that $\Gamma_N\cap \Gamma_C=\emptyset$. The part of the boundary $\Gamma_C$ is in contact with a foundation. We are interested in a mathematical model that describes the evolution of body displacement, stress, and temperature. Let $T>0$ be given, and let $[0,T]$ represent the time interval considered. In what follows, we denote by $x\in\Omega$ and $t\in[0,T]$ the spatial and temporal variables, respectively. We also denote by $u = (u_{i}(x,t))$, $\sigma = (\sigma_{ij}(x,t))$ and $\theta=\theta(x,t)$ the displacement field, the stress field, and the temperature, respectively. These quantities are assumed to satisfy the following constitutive law (for simplicity, we omit the symbols $x$ and $t$ below)
\begin{equation}
\sigma = \mathcal{A}(\varepsilon(\dot{u})) + \mathcal{B}\varepsilon({u}) + C\theta \quad \mbox{in} \ \Omega \times (0,T).\label{classical_sigma}
\end{equation}
 The symbols $\mathcal{A}$, $\mathcal{B}$ and $C$ in (\ref{classical_sigma}) represent the viscosity operator, the elasticity operator, and the heat expansion tensor, respectively. For specific examples of these operators, we refer to Section~\ref{SecNum}. The process is assumed to be dynamic, and hence the equation of motion takes the form of
\begin{equation}
\rho \ddot{u} = \Div \sigma + f_0 \quad \mbox{in} \ \Omega\times (0,T),\label{classical_eq}
\end{equation}
where $\rho$ is the density of the material and $f_0$ is the density of the volume forces. For simplicity of notation, we assume that $\rho = 1$. We assume that the traction force is applied on $\Gamma_N$, that is,
\begin{equation}
\sigma\nu = f_2 \quad \ \mbox{on} \ \Gamma_N\times (0,T).\label{classical_traction}
\end{equation}
The normal compliance condition is described by relation
\begin{equation}
-\sigma_\nu = p(u_\nu) \quad \mbox{on} \ \Gamma_C\times (0,T),\label{classical_bilat}
\end{equation}
where $p_\nu$ is a given non-negative function that vanishes when its argument
is negative.
Moreover, we assume that the frictional law takes the following multivalued, possibly non-monotone form
\begin{equation}
-\sigma_\tau \in \partial j(\dot{u}_\tau) \quad \mbox{on} \ \Gamma_C\times (0,T),\label{classical_frict}
\end{equation}
where $j$ is a locally Lipschitz function, and $\partial j$ stands for its Clarke subgradient; see Section~\ref{Preliminaries}.
In our model, we assume a simplified and linearized version of the energy law that takes the form
\begin{equation}
\rho c_p\dot{\theta} - \div(K\nabla\theta) = c_{ij}\dot u_{i,j} + g \quad \mbox{on} \ \Omega\times (0,T),\label{classical_energy}
\end{equation}
where $c_p$ is the heat capacity, $K$ is the thermal conductivity and $g$ represents the heat sources. We assume that the heat transfer between the body and the foundation is described by the law
\begin{equation}
-K\nabla\theta \cdot \nu = -r(\theta) - h(\|\dot{u}_\tau\|_{\real^d}) \quad \mbox{on} \ \Gamma_C\times (0,T).\label{classical_flux}
\end{equation}
Finally, we impose the initial conditions
\begin{equation}
u(0)=u_0, \ \dot{u}(0)=v_0, \ \theta(0)=\theta_0 \quad \mbox{in} \ \Omega. \label{classical_intial}
\end{equation}
The mechanical problem reads as follows\\

\noindent {\bf Problem ${\cal P_M}$.} {\it Find the displacement $u\colon\Omega\times[0,T]\rightarrow\mathbb{R}^d$, the stress
	$\sigma\colon\Omega\times[0,T]\rightarrow\mathbb{S}^d$ and the temperature $\theta\colon\Omega\times[0,T]\to\real$ that satisfy conditions $(\ref{classical_sigma})$-$(\ref{classical_intial})$.}\\[3mm] 
In the study of Problem ${\cal P_M}$ we need the following assumptions on its data.

\medskip
\noindent
$\underline{H(\mathcal{A})}\ $: The viscosity operator $\mathcal{A}: \Omega \times [0,T]\times \mathbb{S}^d \to \mathbb{S}^d$ satisfies
\lista{
	\item[(a)] $\mathcal{A}(\cdot,\cdot,\varepsilon)$ is continuous on $\Omega \times [0,T]$ for all $\varepsilon \in \mathbb{S}^d$;
	\smallskip
	\item[(b)] $\big(\mathcal{A}(x,t,\varepsilon_1) - \mathcal{A}(x,t,\varepsilon_2)\big):(\varepsilon_1 -\varepsilon_2) \ge m_{\mathcal{A}} \|\varepsilon_1 - \varepsilon_2\|_{\mathbb{S}^d}^2$ for all $\varepsilon_1,\varepsilon_2 \in \mathbb{S}^d$, a.e. $(x,t)\in \Omega\times (0,T)$ with $m_{\mathcal{A}} >0$;
	\smallskip
	\item[(c)] $\|\mathcal{A}(x,t,\varepsilon_1)-\mathcal{A}(x,t,\varepsilon_2)\|_{\mathbb{S}^d}\le L_\mathcal{A} \|\varepsilon_1 - \varepsilon_2\|_{\mathbb{S}^d}$ for all $\varepsilon_1,\varepsilon_2 \in \mathbb{S}^d$, a.e. $(x,t)\in \Omega\times (0,T)$ with $L_\mathcal{A} > 0$.
}

\medskip
\noindent
$\underline{H(\mathcal{B})}\ $: The elasticity operator $\mathcal{B}: \Omega\times \mathbb{S}^d \to \mathbb{S}^d$ is bounded, symmetric, positive fourth-order tensor, i.e.
\lista{
	\item[(a)]  $\mathcal{B}_{ijkl} \in L^\infty(\Omega)$, $1\le i,j,k,l\le d$;
	\smallskip
	\item[(b)] $\mathcal{B}\sigma : \tau = \sigma : \mathcal{B}\tau$ for all $\sigma,\tau\ \in \mathbb{S}^d$, a.e. in $\Omega$;
	\smallskip
	\item[(c)] $\mathcal{B}\tau : \tau \ge 0$ for all $\tau \in \mathbb{S}^d$, a.e. in $\Omega$.
}

\medskip
\noindent
$\underline{H(C)}\ $: The thermal expansion tensor $C = (c_{ij})$ satisfies

\lista{
	\item[(a)] $c_{ij} \in L^\infty (\Omega)$, $1\le i,j\le d$;
	\smallskip
	\item[(b)] $c_{ij}=c_{ji}$,  $1\le i,j\le d$.
}

\medskip
\noindent
$\underline{H(K)}:$ The thermal conductivity operator  
$\displaystyle K \colon \Omega\times [0,T] \times \real^d \to \real^d$ satisfies

\lista{
	\item[(a)]
	$K(\cdot,\cdot, \xi)$ is continuous on $\Omega \times [0,T]$ for all $\xi \in \real^d$;
	\smallskip
	\item[(b)]
	$\displaystyle 
	\big( K(x, t, \xi_1) - K (x, t, \xi_2) \big) \cdot (\xi_1 - \xi_2) \ge 
	m_K \, \| \xi_1 - \xi_2 \|_{\real^d}^2$
	for all $\xi_1$, $\xi_2 \in \real^d$, a.e. $(x, t) \in \Omega\times(0,T) $ with $m_K > 0$;  
	\smallskip
	\item[(c)] $\|K(x,t,\xi_1)-K(x,t,\xi_2)\|_{\real^d}\le L_K \|\xi_1-\xi_2\|_{\real^d}$ for all $\xi_1,\xi_2\in \real^d$, a.e. $(x,t)\in \Omega\times (0,T)$ with $L_K>0$.
}

\medskip

\noindent
$\underline{H(p)}:$ The function  
$\displaystyle p \colon \Gamma_C\times [0,T]\times\real \to \real$ satisfies 
\lista{
	\item[(a)]
	$p(\cdot,\cdot,\xi)$ is measurable on $\Gamma_C\times (0,T)$ for all 
	$\xi \in \mathbb{R}^d$;
	\smallskip
	\item[(b)]
	$\displaystyle 
	| p (x, t, \xi_1) - p (x, t, \xi_2) | \le L_p \, | \xi_1 - \xi_2 |$
	for all $\xi_1$, $\xi_2 \in \real$, a.e. $(x, t) \in \Gamma_C\times (0,T)$ with $L_p > 0$.
}
\medskip

\noindent
$\underline{H(j)}\  $: The potential  
$j\colon \Gamma_C\times [0,T] \times \mathbb{R}^d \to \real$ satisfies


\lista{
	\item[(a)]
	$j(\cdot, \cdot, \xi)$ is measurable on $\Gamma_C\times (0,T)$ for all 
	$\xi \in \mathbb{R}^d$;
	\smallskip
	\item[(b)]
	$j(x, t,\cdot)$ is locally Lipschitz on $\real^d$ for a.e. $(x, t) \in \Gamma_C\times (0,T)$;
	\smallskip
	\item[(c)]
	$\| \partial j(x, t, \xi) \|_{\mathbb{R}^d} \le c_j $ for all $\xi \in \mathbb{R}^d$, 
	a.e. $(x, t) \in \Gamma_C\times (0,T)$ with  $c_j \ge 0$;
	\smallskip
	\item[(d)]
	$(\zeta_1 - \zeta_2) \cdot (\xi_1 - \xi_2) \ge -
	m_j \| \xi_1 - \xi_2 \|^2_{\mathbb{R}^d}$ for all 
	$\zeta_i \in \partial j(x, t, \xi_i)$, $\xi_i \in \mathbb{R}^d$, $i =1$, $2$, 
	a.e. $(x, t) \in \Gamma_C\times (0,T)$ with $m_j \ge 0$.
}

\medskip

\noindent
$\underline{H(r)}:$ The function  
$\displaystyle r \colon \Gamma_C\times\real \to \real$ satisfies 

\lista{
	\item[(a)]
	$r(\cdot,\zeta) \in L^2(\Gamma_C)$ for all $\zeta \in \real$;
	\smallskip
	\item[(b)]
	$\displaystyle 
	| r (x, \zeta_1) - r (x, \zeta_2) | \le L_r \, | \zeta_1 - \zeta_2 |$
	for all $\zeta_1$, $\zeta_2 \in \real$, a.e. $x \in \Gamma_C$ with $L_r > 0$.
}

\medskip

\noindent
$\underline{H(h)}:$ The function 
$\displaystyle h \colon \Gamma_C \times [0,\infty) \to [0,\infty)$ satisfies 

\lista{
	\item[(a)]
	$h(\cdot, \zeta) \in L^2(\Gamma_C)$ for all $\zeta \in [0,\infty)$;
	\smallskip
	\item[(b)]
	$\displaystyle 
	| h(x, \zeta_1) - h(x, \zeta_2) | \le L_h \, | \zeta_1 - \zeta_2 |$
	for all $\zeta_1$, $\zeta_2 \in [0,\infty)$, a.e. $x \in \Gamma_C$ with $L_h > 0$.
}

\medskip

\noindent $\underline{H(f)}\ $: The volume forces and traction densities satisfy
\[
f_0 \in C(0,T;L^2(\Omega;\real^d)), \quad f_2 \in C(0,T;L^2(\Gamma_N;\real^d)).
\]

\medskip

\noindent $\underline{H(g)}\ $: The heat source satisfies
\[
g \in C(0,T;L^2(\Omega)).
\]

\medskip

\noindent We remark that the result presented in the sequel can be obtained under a weaker assumption than $H(j)(c)$. Namely it is enough to assume a sublinear growth of the subgradient $\partial j$ with respect to the last variable. However, for the sake of simplicity we assume that it is bounded.\\ 

\noindent In order to provide a variational formulation of the problem, we define spaces 
\begin{equation}
E=H^1(\Omega), \quad V=H^1(\Omega;\real^d),\nonumber
\end{equation}
where $E$ is equipped with the classical $H^1$ norm, and in $V$ we define the inner product
\[
(u,v)_V = (u,v)_{H} + (\varepsilon(u),\varepsilon(v))_{Q} \quad \mbox{for\,\,all} \ u,v\in V
\]
and the corresponding norm $\|u\|^2:=(u,u)_V$ for all $u\in V$. It follows that the norms $\|\cdot\|_V$ and $\|\cdot\|_{H^1(\Omega;\real^d)}$ are equivalent. We know that the embedding  $i\colon E \to H^{1-\delta}(\Omega)$ is compact for $\delta\in (0,\frac{1}{2})$, and the trace operator  $\bar{\gamma}\colon H^{1-\delta}(\Omega) \to L^2(\Gamma)$ is continuous. Hence, we see that $\gamma=\bar{\gamma}i\colon E \to L^2(\Gamma)$ is compact. We omit the notation of the embedding and the trace operator, and write simply $v$ instead of $\gamma v$ for $v\in E$. We proceed similarly for the space $V$. We introduce the following spaces of vector-valued functions
\[
\mathbb{V} = \{v\in L^2(0,T;V)\mid \dot v\in L^2(0,T;V^*) \}, \quad \mathbb{E} = \{v\in L^2(0,T;E)\mid \dot v\in L^2(0,T;E^*) \},
\]
where $\dot v$ denotes time derivative of $v$ in the sense of distributions.

Next, we define operators $A\colon(0,T)\times V \to V^*$, $B\colon V\to V^*$, $C_1\colon E\to V^*$, $C_2\colon (0,T)\times E\to E^*$, $C_3\colon V\to E^*$ by
\begin{align*}
& 
\langle A(t,u), v \rangle_V  =({\cal A}(t,\varepsilon(u)),\varepsilon(v))_{Q} \quad {\rm for\,\,all} \ u,v\in V, \ t\in(0,T), \\[2mm]
& 
\langle Bu,v \rangle_V = ({\cal B}\varepsilon(u),\varepsilon(v))_{Q}  \quad {\rm for\,\,all} \ u,v\in V, \\ 
& 
\langle C_1 \theta, v \rangle_V = \int_{\Omega} c_{ij} v_{i,j} \theta \,dx \quad {\rm for\,\,all} \ v\in V, \theta \in E, \\ 
& 
\langle C_2 (t,\theta), \eta \rangle_E = ( K(t,\nabla \theta),\nabla \eta )_{L^2(\Omega)} \quad {\rm for\,\,all} \ \theta,\eta\in E,\\ 
&
\langle C_3 v,\eta \rangle_E = -\int_{\Omega} c_{ij} v_{i,j} \eta \,dx \quad  {\rm for\,\,all} \ v\in V, \eta \in E. 
\end{align*}
We define the functions $f \colon (0, T) \to V^*$ and $\tilde{g}\colon (0,T)\to E^*$ for all $v\in V$, $\eta\in E$, a.e. $t\in (0,T)$ by
\begin{align*}
&\langle f(t), v \rangle_{V^* \times V} = \langle f_0(t), v \rangle_{H} + 
\langle f_2(t), v \rangle_{L^2(\Gamma_N; \real^d)},\\[2mm]
&\dual{\tilde{g}(t),\eta}{E}=\skalar{g(t),\eta}_{L^2(\Omega)}.
\end{align*}

\noindent We also introduce an additional hypothesis on the data

\medskip

\noindent
$\underline{H_0}\ :$ $u_0\in V$, $v_0\in H$, $\theta_0\in E$.

\medskip

\indent Now we are in a position to present a weak formulation of Problem ${\cal P_M}$.

\newpage

\noindent {\bf Problem ${\cal P_E}$.} {\it 
	Find a displacement $u\in L^2(0,T;V)$ with $\dot u\in \mathbb{V}$, and a temperature $\theta \in \mathbb{E}$ such that
\begin{align}
&
\label{eqKB_39}\langle \ddot{u}(t) + A(t,\dot{u}(t))+ Bu(t) + C_1\theta(t)-f(t),w\rangle_V =  - \int_{\Gamma_C} p(u_\nu) w_\nu \, d\Gamma\\
&\nonumber\hspace{2.3cm}
 - \int_{\Gamma_C} \xi(t) \cdot w_\tau \, d\Gamma \qquad \mbox{for all} \ w\in V, \ \mbox{a.e.} \ t\in (0,T),\\
&
\label{eqKB_40}\langle \dot{\theta}(t) + C_2(t,\theta(t)) + C_3 \dot{u}(t) -\tilde{g}(t),\eta \rangle_E = \int_{\Gamma_C} r(\theta)\eta \,d\Gamma \\
&\nonumber\hspace{2.3cm}
+\int_{\Gamma_C}h(\|\dot{u}_\tau(t)\|_{\real^d})\eta\,d\Gamma\quad \mbox{for all} \ \eta\in E,  \ \mbox{a.e.} \ t\in (0,T),\\
&
\xi\in \partial j(\dot{u}_\tau),\ \mbox{a.e.} \ (x,t)\in\Gamma_C\times (0,T),\label{eqKB_14}\\
&
\nonumber u(0)=u_0, \ \dot{u}(0)=v_0, \ \theta(0)=\theta_0.
\end{align}
}	
\hspace{-.17cm}Problem ${\cal P_E}$ can be obtained from equations (\ref{classical_eq}) and (\ref{classical_energy}) by multiplying them by test functions $w\in V$ and $\eta\in E$, respectively, performing the Green formulas (\ref{Green_a}) and (\ref{Green}), applying conditions \eqref{classical_sigma}--\eqref{classical_intial}, using definitions of operators above and denoting $\xi=\sigma_\tau$. Every solution of Problem ${\cal P_E}$ is called a weak solution of Problem ${\cal P_M}$.

	Existence and uniqueness of the solution to Problem ${\cal P_E}$ can be proved following the proof of Theorem 9 in \cite{SZAF}
	 and using Corollary~\ref{Erhilng_Gamma}. Since this paper is devoted to numerical analysis, we omit this proof. For more details on non-clamped problems see e.g. \cite{Chau2015} and for problems with no smallness assumption see e.g. \cite{BARTOSZ2021105940}.

\begin{Corollary}\label{Erhilng_Gamma}
	Lemma~\ref{Erhling} holds for spaces $X=H^1(\Omega)$, $Y=H^{1-\delta}(\Omega)$, $Z=L^2(\Omega)$ for $\delta \in (0,1/2)$. Moreover, from the trace theorem, we have $\|v\|_{L^2(\Gamma)}\le C\|v\|_{H^{1-\delta}(\Omega)}$, hence for all $\varepsilon>0$ there exists $C(\varepsilon)>0$ such that 
	\[
	\|v\|_{L^2(\Gamma)} \le \varepsilon \|v\|_E + C(\varepsilon)\|v\|_{L^2(\Omega)} \quad \mbox{for all} \ v\in E.
	\]
	Similarly, we obtain
	\[
	\|v\|_{L^2(\Gamma;\real^d)} \le \varepsilon \|v\|_V + C(\varepsilon)\|v\|_{L^2(\Omega;\real^d)} \quad \mbox{for all} \ v\in V.
	\]

\end{Corollary}

\section{Discretization and error estimates}\label{SecError}
	
	In this section we consider a fully discrete approximation of Problem ${\cal P_E}$ and examine error estimates for the approximate solution and the exact one. 
	
	Let $E^h$ and $V^h$ be finite-dimensional subspaces of $E$ and $V$, respectively, with $h>0$ denoting the spatial discretization parameter. 
	Let $u_0^h,v_0^h\in V^h$ and $\theta_0^h\in E^h$ be suitable approximations of $u_0,v_0,\theta_0$. We define a uniform partition of $[0,T]$ denoted by $0=t_0<t_1<\ldots<t_N=T$, where $N\in\mathbb{N}$, $t_j=jk$ for $j=1,...,N$ and  $k=\frac{T}{N}$ is the time step. For a continuous function $\varphi$ we use the notation 
$\varphi_n=\varphi(t_n)$. For a sequence $\{z_n\}_{n=0}^N$ we denote by $\delta z_n=(z_n-z_{n-1})/k$ for $n=1,\ldots,N$ its backward divided difference.\\
Let $(u,\theta)$ be a solution of Problem ${\cal P_E}$ and  $v=\dot u$. The fully discrete approximation of Problem ${\cal P_E}$ is the following.\\[3mm]
\noindent {\bf Problem ${\cal P}_{\cal E}^{hk}$.} {\it
	Find sequences of velocity $\{v_n^{hk}\}_{n=1}^N\subset V^h$ and temperature $\{\theta_n^{hk}\}_{n=1}^N \subset E^h$ such that there exists $\{\xi_n^{hk}\}_{n=1}^N \subset L^2(\Gamma_C;\real^d)$ satisfying for $j=1,\ldots,N$
	\begin{align}
	&\label{eqKB_19}
	\hspace{-.2cm}\langle \delta v_j^{hk} + A_j(v_j^{hk}) + Bu_j^{hk} + C_1 \theta_j^{hk} -f_j, w^h\rangle_V=  - \int_{\Gamma_C} p(u_{j\nu}^{hk})\cdot w_\nu^h \, d\Gamma\\
    &\nonumber\hspace{2.3cm}
     - \int_{\Gamma_C} \xi_j^{hk} \cdot w_\tau^h \, d\Gamma 
     \qquad \mbox{for all} \ w^h\in V^h, \ \mbox{a.e.} \ t\in (0,T),\\[2mm]
	&
	\hspace{-.2cm}\langle \delta\theta_j^{hk} + C_2(\theta_j^{hk}) + C_3 v_j^{hk}-\tilde{g}_j,\eta^h\rangle_E = \int_{\Gamma_C} r(\theta_j^{hk}) \, \eta^h\,d\Gamma\label{eqKB_20}\\[2mm]
	&\hspace{2.3cm} +\int_{\Gamma_C}h(\|v_{j\tau}^{hk}\|_{\real^d})\eta^h\,d\Gamma \qquad \mbox{for all} \ \eta^h\in E^h,\nonumber \\[2mm] 
	&
	\hspace{-.2cm}\xi_j^{hk}\in \partial j(v_{j\tau}^{hk}),\label{eqKB_21} 
	\end{align}
	where $v_0^{hk}=v_0^h, \theta_0^{hk} = \theta_0^h$ and $\{u_n^{hk}\}_{n=1}^N \subset V^h$ is given by
	\begin{equation}\label{unhk}
	u_n^{hk} = u_0^{hk} +\sum_{j=1}^n kv_j^{hk} \quad\text{for}\,\,n=1,...,N
	\end{equation}
	with $u_0^{hk}=u_0^h$.
}

\medskip

Now we deal with the error analysis for the discrete scheme, Problem ${\cal P}_{\cal E}^{hk}$.\\[2mm]
Let $\{u_n^{hk}\}_{n=1}^N$, $\{v_n^{hk}\}_{n=1}^N$ and $\{\theta_n^{hk}\}_{n=1}^N$ denote the discrete approximation of displacement $u$, velocity $v$ and temperature $\theta$, respectively, derived by the numerical scheme in Problem ${\cal P}_{\cal E}^{hk}$. We prove the following estimate.

\begin{Theorem}\label{main2}
	Let the assumptions  $H({\cal A})$, $H({\cal B})$, $H(C)$, $H(K)$, $H(j)$, $H(r)$, $H(h)$, $H(f)$, $H(g)$, and $H_0$ hold. Moreover, we impose the following  additional regularity conditions  
	$u\in C^2(0,T;H)\cap C^1(0,T;V)$, $\dot{u}_\tau \in C(0,T;L^2(\Gamma_C;\real^d))$, $\theta\in C(0,T;L^2(\Gamma_C))$. Then, there exists $C>0$ such that for all $\{w_n^h\}_{n=1}^N\subset V^h$ and $\{\eta_n^h\}_{n=1}^N\subset E^h$ we have
	\begin{align}\label{KB2}
	&
	\max_{1\le n\le N} \{\|v_n-v_n^{hk}\|_H^2 + \|\theta_n-\theta_n^{hk}\|_{L^2(\Omega)}^2\} + k\sum_{j=1}^N \|v_j-v_j^{hk}\|_V^2 + k\sum_{j=1}^N \|\theta_j-\theta_j^{hk}\|_E^2\nonumber   \\
	&
	\le C\bigg[ k\sum_{j=1}^N (\|\dot{v}_j -\delta v_j \|_H^2 + \|v_j-w_j^h\|_V^2) + \max_{1\le n\le N} \|v_{n\tau} - w_{n\tau}^h\|_{L^2(\Gamma_C;\real^d)} \nonumber \\
	&
	+ \frac{1}{k}\sum_{j=1}^{N-1} \|(v_j-w_j^h) - (v_{j+1} - w_{j+1}^h )\|_H^2 + \max_{1\le n\le N} \|v_n-w_n^h\|_H^2\nonumber \\
	&
	+ \|u_0-u_0^h\|_V^2 + \|v_0-v_0^h\|_H^2 + \|\theta_0-\theta_0^h\|^2_{L^2(\Omega)}+ k^2\|u\|_{H^2(0,T;V)} \nonumber \\
	&
	+k\sum_{j=1}^N (\|\dot{\theta}_j-\delta \theta_j\|_{L^2(\Omega)   } +\|\theta_j - \eta_j^h\|_E^2)  \nonumber \\
	&
	+\frac{1}{k}\sum_{j=1}^{N-1} \|(\theta_j -\eta_j^h) - (\theta_{j+1} -\eta_{j+1}^h)\|_{L^2(\Omega)}^2 + \max_{1\le n\le N} \|\theta_n-\eta_n^h\|^2_{L^2(\Omega)}\bigg].
	\end{align}
\end{Theorem}

\noindent {\bf Proof.} Take $w^h\in V^h$ and $\eta^h\in E^h$ as test functions in both Problem ${\cal P_E}$ and Problem ${\cal P}_{\cal E}^{hk}$, respectively, and subtract relevant equations to obtain 
\begin{align}
&\nonumber
( \dot{v}_j-\delta v_j^{hk}, w^h)_H + \langle A_j(v_j)-A_j(v_j^{hk}),w^h\rangle_V + \langle B(u_j-u_j^{hk}),w^h\rangle_V  \\[2mm]
&
\nonumber+\langle C_1 (\theta_j- \theta_j^{hk}), w^h\rangle_V + \int_{\Gamma_C} \big(p(u_{j\nu}) -p(u_{j\nu}^{hk}) \big) \cdot w_\nu ^h\, d\Gamma + \int_{\Gamma_C} (\xi_j -\xi_j^{hk} ) \cdot w_\tau ^h\, d\Gamma \\[2mm]
&\nonumber
+( \dot{\theta}_j-\delta \theta_j^{hk},\eta^h )_{L^2(\Omega)} + \langle C_2\theta_j - C_2 \theta_j^{hk},\eta^h\rangle_E+ \langle C_3(v_j - v_j^{hk}),\eta^h\rangle_E   \\[2mm]
&
+\int_{\Gamma_C} \big(r(\theta_j^{hk}) -r(\theta_j)\big)\, \eta^h\, d\Gamma+\int_{\Gamma_C}\left(h(\|v_{j\tau}^{hk}\|_{\real^d})-h(\|v_{j\tau}\|_{\real^d})\right)\eta^h\,d\Gamma = 0 \label{eq_1_zero}
\end{align}
	
\noindent for $j=1,\ldots,N.$ Let $\{w_j^h\}_{j=1}^N\subset V$ and $\{\eta_j^h\}_{j=1}^N\subset E$ be arbitrary. Applying \eqref{eq_1_zero} with $w^h = v_j^{hk}$, $\eta^h=\theta_j^{hk}$ and then with  $w^h=w_j^h$, $\eta^h=\eta_j^h$, we find
\begin{align*}
&
( \dot{v}_j - \delta v_j^{hk},v_j- v_j^{hk})_H + \langle A_j(v_j) - A_j(v_j^{hk}),v_j-v_j^{hk}\rangle_V \\[2mm]
&
+ \langle B(u_j - u_j^{hk}),v_j-v_j^{hk}\rangle_V + \langle C_1(\theta_j - \theta_j^{hk}),v_j-v_j^{hk}\rangle_V + \langle \dot{\theta}_j-\delta \theta_j^{hk},\theta_j-\theta_j^{hk}\rangle_E \\[2mm]
&
+ \langle C_2 \theta_j - C_2\theta_j^{hk},\theta_j-\theta_j^{hk}\rangle_E + \langle C_3(v_j-v_j^{hk}),\theta_j-\theta_j^{hk}\rangle_E  \\[2mm]
&
+ \int_{\Gamma_C} \big(p(u_{j\nu}) - p(u_{j\nu}^{hk}) \big)\cdot (v_{j\nu}-v_{j\nu}^{hk})\,d\Gamma + \int_{\Gamma_C} (\xi_j - \xi_j^{hk})\cdot (v_{j\tau}-v_{j\tau}^{hk})\,d\Gamma \\[2mm]
&
+\int_{\Gamma_C} \big(r(\theta_j^{hk}) - r(\theta_j)\big)\, (\theta_j-\theta_j^{hk})\,d\Gamma +\int_{\Gamma_C}\left(h(\|v_{j\tau}^{hk}\|_{\real^d})-h(\|v_{j\tau}\|_{\real^d})\right)(\theta_j-\theta_j^{hk})\,d\Gamma \\
&
=( \dot{v}_j - \delta v_j^{hk},v_j- w_j^{h})_H + \langle A_j(v_j) - A_j(v_j^{hk}),v_j-w_j^h\rangle_V \\[2mm]
&
+ \langle B(u_j - u_j^{hk}),v_j-w_j^h\rangle_V + \langle C_1(\theta_j - \theta_j^{hk}),v_j-w_j^h\rangle_V + \langle \dot{\theta}_j-\delta \theta_j^{hk},\theta_j-\eta_j^h\rangle_E \\[2mm]
&
+ \langle C_2 \theta_j - C_2\theta_j^{hk},\theta_j-\eta_j^h\rangle_E + \langle C_3(v_j-v_j^{hk}),\theta_j-\eta_j^h\rangle_E  \\[2mm]
&+ \int_{\Gamma_C} \big(p(u_{j\nu}) - p(u_{j\nu}^{hk}) \big)\cdot (v_{j\nu}-w_{j\nu}^{hk})\,d\Gamma + \int_{\Gamma_C} (\xi_j - \xi_j^{hk}) \cdot(v_{j\tau}-w_{j\tau}^{h})\,d\Gamma\\[2mm]
&+ \int_{\Gamma_C} \big(r(\theta_j^{hk}) - r(\theta_j)\big)\, (\theta_j-\eta_j^{h})\,d\Gamma +\int_{\Gamma_C}\left(h(\|v_{j\tau}^{hk}\|_{\real^d})-h(\|v_{j\tau}\|_{\real^d})\right)(\theta_j-\eta_j^{h})\,d\Gamma  . 
\end{align*}

\noindent
Using the fact that $\langle C_1 \eta,w\rangle_V +\langle C_3w,\eta\rangle_E = 0$ for all $w\in V$, $\eta\in E$, we get
\begin{align}
&\nonumber
( \delta v_j -\delta v_j^{hk},v_j-v_j^{hk} )_H + \langle A_j (v_j)-A_j (v_j^{hk}),v_j-v_j^{hk}\rangle_V\\[2mm]
&\nonumber + ( \delta \theta_j-\delta\theta_j^{hk},\theta_j -\theta_j^{hk} )_{L^2(\Omega)} + \langle C_2\theta_j -C_2\theta_j^{hk},\theta_j - \theta_j^{hk}\rangle_E \\[2mm]
&\nonumber
+ \int_{\Gamma_C} (\xi_j - \xi_j^{hk})\cdot (v_{j\tau}-v_{j\tau}^{hk})\,d\Gamma + \int_{\Gamma_C} \big(r(\theta_j^{hk}) - r(\theta_j)\big) (\theta_j-\theta_j^{hk})\,d\Gamma \\[2mm]
&\nonumber
+ \int_{\Gamma_C}\left(h(\|v_j^{hk}\|_{\real^d})-h(\|v_j\|_{\real^d})\right)(\theta_j-\theta_j^{hk})\,d\Gamma \\[2mm]
&\nonumber
=( \delta v_j- \delta v_j^{hk},v_j- w_j^h )_H + ( \dot{v}_j -\delta v_j , (v_j -w_j^h) + (v_j^{hk}  -v_j) )_H \\[2mm]
&\nonumber
+ \langle A_j (v_j)-A_j (v_j^{hk}),v_j - w_j^h\rangle_V + \langle B(u_j-u_j^{hk}),(v_j-w_j^h) + (v_j^{hk} -v_j)\rangle_V  \\[2mm]
&\nonumber
+\langle C_1(\theta_j - \theta_j^{hk}),v_j-w_j^h\rangle_V + ( \delta \theta_j -\delta \theta_j^{hk},\theta_j - \eta_j^h )_{L^2(\Omega)}  \\[2mm]
&\nonumber
+( \dot{\theta}_j-\delta \theta_j,(\theta_j -\eta_j^h) + (\theta_j^{hk}-\theta_j) )_{L^2(\Omega)}+ \langle C_2\theta_j - C_2 \theta_j^{hk},\theta_j-\eta_j^h\rangle_E 
\\[2mm]
&
+ \langle C_3(v_j-v_j^{hk}),\theta_j - \eta_j^h \rangle_E +   \int_{\Gamma_C}\big(p(u_{j\nu}) - p(u_{j\nu}^{hk}) \big)\cdot \big((v_{j\nu}-w_{j\nu}^{hk}) + (v_{j\nu}^{hk} - v_{j\nu})\big)\,d\Gamma \nonumber\\[2mm]
&\nonumber + \int_{\Gamma_C} (\xi_j - \xi_j^{hk})\cdot (v_{j\tau}-w_{j\tau}^{h})\,d\Gamma + \int_{\Gamma_C} \big(r(\theta_j^{hk}) - r(\theta_j)\big) (\theta_j-\eta_j^{h})\,d\Gamma\\[2mm]
&+\int_{\Gamma_C}\left(h(\|v_{j\tau}^{hk}\|_{\real^d})-h(\|v_{j\tau}\|_{\real^d})\right)(\theta_j-\eta_j^{h})\,d\Gamma. \label{eqKB_36}
\end{align}

\noindent From the formula $2\skalar{a-b,a}_H = \|a-b\|_H^2+\|a\|_H^2-\|b\|_H^2$ for $a,b\in H$, we get
\begin{align}
&\label{eqKB_37}
(\delta v_j -\delta v_j^{hk},v_j-v_j^{hk} )_H \geq\frac{1}{2k}(\|v_j-v_j^{hk}\|_H^2 -\|v_{j-1} -v_{j-1}^{hk}\|_H^2)   , \\[2mm]
&
( \delta \theta_j -\delta \theta_j^{hk},\theta_j-\theta_j^{hk} )_{L^2(\Omega)}\geq \frac{1}{2k}(\|\theta_j-\theta_j^{hk}\|_{|L^2(\Omega)}^2 -\|\theta_{j-1} -\theta_{j-1}^{hk}\|_{L^2(\Omega)}^2).
\end{align}

\noindent
By $H({\cal A})(b)$, $H(j)(d)$ and Corollary~\ref{Erhilng_Gamma},  we get
\begin{align}
&\nonumber
\langle A_j (v_j) - A_j (v_j^{hk}),v_j-v_j^{hk}\rangle_V + \int_{\Gamma_C} (\xi_j-\xi_j^{hk})\cdot (v_{j\tau}-v_{j\tau}^{hk})\, d\Gamma  \\[2mm]
&\nonumber
\ge m_{\mathcal{A}}\|v_j-v_j^{hk}\|_V^2 - m_{\mathcal{A}}\|v_j-v_j^{hk}\|_{H}^2 - m_j\|v_j-v_j^{hk}\|_{L^2(\Gamma_C;\real^d)}^2\\[2mm]
&
\ge (m_{\mathcal{A}} - \varepsilon)\|v_j-v_j^{hk}\|_V^2 - \big(m_{\mathcal{A}} +  C(\varepsilon)\big)\|v_j-v_j^{hk}\|_{H}^2.
\end{align}
Similarly, we obtain
\begin{align}
&\nonumber
\langle C_2\theta_j - C_2\theta_j^{hk},\theta_j-\theta_j^{hk}\rangle_E + \int_{\Gamma_C} \big(r(\theta_j^{hk}) -r(\theta_j)\big)(\theta_j-\theta_j^{hk})\, d\Gamma \\[2mm]
&\nonumber
+ \int_{\Gamma_C} \big(h(\|v_j^{hk}\|_{\mathbb{R}^d})-h(\|v_j\|_{\mathbb{R}^d})\big)(\theta_j-\theta_j^{hk})\, d\Gamma  \\[2mm]
&\nonumber
\ge m_K\|\theta_j-\theta_j^{hk}\|_E^2 - m_K\|\theta_j-\theta_j^{hk}\|_{L^2(\Omega)}^2\\[2mm]
&\nonumber - (L_r+\frac{1}{2})\|\theta_j-\theta_j^{hk}\|_{L^2(\Gamma)}^2 - \frac{L_h^2}{2}\|v_j-v_j^{hk}\|_{L^2(\Gamma_C)}^2 \\[2mm]
&\nonumber
\ge (m_K-\varepsilon )\|\theta_j-\theta_j^{hk}\|_E^2 -\big(m_K+C(\varepsilon)\big) \|\theta_j-\theta_j^{hk}\|_{L^2(\Omega)}^2\\[2mm] 
&-\varepsilon \|v_j-v_j^{hk}\|_V^2 
 -  C(\varepsilon)\|v_j-v_j^{hk}\|_H^2.
\end{align}

 \noindent
Next, from the Lipschitz continuity of ${\cal A}$ and $K$, we have
\begin{equation}
\langle A_j (v_j) - A_j (v_j^{hk}),v_j - w^h_j\rangle_V \le   \varepsilon\|v_j-v_j^{hk}\|_V^2+C(\varepsilon) \|v_j -w^h_j\|^2_V,
\end{equation}
\begin{equation}
\langle C_2\theta_j -C_2\theta_j^{hk},\theta_j-\eta^h_j\rangle_E \le  \varepsilon\|\theta_j -\theta_j^{hk}\|_E^2+C(\varepsilon) \|\theta_j-\eta^h_j\|_E^2.
\end{equation}
By the linearity of $B$, $C_1$ and $C_3$, we obtain
\begin{align}
&\nonumber\langle B(u_j-u_j^{hk}),(v_j-w_j^h) + (v_j^{hk} -v_j)\rangle_V\leq\varepsilon\|v_j-v_j^{hk}\|_V^2\\[2mm]
&+C(\varepsilon)\|u_j-u_j^{hk}\|^2_V+C\|v_j-w_j^h\|^2_V,\\[2mm]
&\langle C_1(\theta_j - \theta_j^{hk}),v_j-w_j^{h}\rangle_V \le  \varepsilon\|\theta_j - \theta_j^{hk}\|_E^2+C(\varepsilon) \|v_j-w_j^{h}\|_V^2,\\[2mm]
&\langle C_3 (v_j-v_j^{hk}),\theta_j-\eta_j^{h}\rangle_E \le  \varepsilon\|v_j-v_j^{hk}\|_V^2+C(\varepsilon) \|\theta_j-\eta_j^{h}\|_E^2.
\end{align}
From $H(p)$, and the continuity of the trace operator $\gamma\colon V\to L^2(\Gamma;\real^d)$, we deduce
\begin{align}
&\nonumber\int_{\Gamma_C} \big(p(u_{j\nu}) - p(u_{j\nu}^{hk}) \big)\cdot \big((v_{j\nu}-w_{j\nu}^h) + (v_{j\nu}^{hk} -v_{j\nu})\big)\,d\Gamma \\[2mm]
&\qquad \leq \varepsilon \|v_{j\nu}-v_{j\nu}^{hk}\|_{L^2(\Gamma_C)}^2 + C(\varepsilon)\|u_{j\nu}-u_{j\nu}^{hk}\|_{L^2(\Gamma_C)}^2+C\|v_{j\nu}-w_{j\nu}^h\|_{L^2(\Gamma_C)}^2 \\[2mm]
&\qquad \leq \varepsilon \|v_j-v_j^{hk}\|_V^2 + C(\varepsilon)\|u_j-u_j^{hk}\|_V^2+C\|v_j-w_j^h\|_V^2.
\end{align}
By $H(j) (c)$, $H(r) (b)$, $H(h) (b)$, and the continuity of the trace operator $\gamma\colon E\to L^2(\Gamma)$, we get
\begin{align}
&\int_{\Gamma_C} (\xi_j^{hk} - \xi_j)\cdot(v_{j\tau} - w_{j\tau}^h)\,d\Gamma \le 2c_j\sqrt{{\rm meas}(\Gamma_C)} \|v_{j\tau}-w_{j\tau}^h\|_{L^2(\Gamma_C;\real^d)}\nonumber\\[2mm]
& \quad\leq C\|v_j-w_j^h\|_{L^2(\Gamma_C;\real^d)},\label{eqKB_38a}\\[2mm]
&\nonumber\int_{\Gamma_C} \big(r(\theta_j^{hk}) - (r(\theta_j)\big)(\theta_j - \eta_j^h)\,d\Gamma\le L_r\|\theta_j^{hk}-\theta_j\|_{L^2(\Gamma)} \|\theta_j-\eta_j^h\|_{L^2(\Gamma_C)}\\[2mm]
&\quad\leq\varepsilon\|\theta_j^{hk}-\theta_j\|_E^2+C(\varepsilon)\|\theta_j-\eta_j^h\|_E^2,\\[2mm]
&\nonumber\int_{\Gamma_C}\left(h(\|v_{j\tau}^{hk}\|_{\real^d})-h(\|v_{j\tau}\|_{\real^d})\right)(\theta_j-\eta_j^{h})\,d\Gamma\\[2mm]
&\quad\leq\varepsilon\|v_j^{hk}-v_j\|_V^2+C(\varepsilon)\|\theta_j-\eta_j^h\|^2_{E}.\label{eqKB_38}
\end{align}
The symbol ${\rm meas}(\Gamma_C)$ in (\ref{eqKB_38a}) represents the Lebesgue surface measure of $\Gamma_C$.
Finally, we estimate
\begin{align}
&( \dot{v}_j -\delta v_j , (v_j -w_j^h) + (v_j^{hk}  -v_j) )_H\nonumber\\[2mm]
&\leq\|\dot{v}_j -\delta v_j\|^2_H+\frac{1}{2}\|v_j -w_j^h\|^2_V+\frac{1}{2}\|v_j^{hk}  -v_j\|^2_H,\label{eqKB_38b}\\[2mm]
&( \dot{\theta}_j-\delta \theta_j,(\theta_j -\eta_j^h) + (\theta_j^{hk}-\theta_j) )_{L^2(\Omega)}\nonumber\\[2mm]
&\leq\|\dot{\theta}_j-\delta \theta_j\|^2_{L^2(\Omega)}+\frac{1}{2}\|\theta_j -\eta_j^h\|^2_E+\frac{1}{2}\|\theta_j^{hk}-\theta_j\|_{L^2(\Omega)}^2.\label{eqKB_38c} 
\end{align}
Applying (\ref{eqKB_37})-(\ref{eqKB_38c}) in (\ref{eqKB_36}), we deduce
\begin{align*}
&
\frac{1}{2k} (\|v_j-v_j^{hk}\|_H^2 -\|v_{j-1}-v_{j-1}^{hk}\|_H^2) + \left(m_{\mathcal{A}}-7\varepsilon\right)\|v_j-v_j^{hk}\|_V^2\\[2mm]
&
+ \frac{1}{2k}(\|\theta_j-\theta_j^{hk}\|_{L^2(\Omega)}^2- \|\theta_{j-1} -\theta_{j-1}^{hk}\|_{L^2(\Omega)}^2) + \left(m_K-4\varepsilon\right)\|\theta_j-\theta_j^{hk}\|_E^2 \\[2mm]
&
\le C \left(\|\dot{v}_j -\delta v_j\|_H^2 + \|v_j-w_j^h\|_V^2 + \|u_j-u_j^{hk}\|_V^2 + \|\dot{\theta}_j -\delta \theta_j\|_{L^2(\Omega)}^2 
\right.\\[2mm]
&\left.+\|\theta_j-\eta_j^h\|_E^2+\|v_j -w_j^h\|_{L^2(\Gamma_C;\real^d)}
 + \|v_j-v_j^{hk}\|_{H}^2 \right.\\[2mm]
&\left.+\|\theta_j-\theta_j^{hk}\|_{L^2(\Omega)}^2 \right)  
+( \delta v_j-\delta v_j^{hk},v_j-w_j^h)_H
+ ( \delta\theta_j -\delta\theta_j^{hk},\theta_j -\eta_j^h)_{L^2(\Omega)}.
\end{align*}

\noindent
Summing up the last inequality from $j=1$ to $j=n$, and taking into account that $v_0=\dot{u}(0)=v_0$, we obtain the following
\begin{align*}
&
\|v_n-v_n^{hk}\|_H^2 + 2k(m_{\mathcal{A}}-7\varepsilon) \sum_{j=1}^n \|v_j-v_j^{hk}\|_V^2 + \|\theta_n-\theta_n^{hk}\|_{L^2(\Omega)}^2 \\
&
+ 2k(m_K-4\varepsilon) \sum_{j=1}^n \|\theta_j-\theta_j^{hk}\|_E^2 \le \|v_0-v_0^{hk}\|_V^2 + \|\theta_0-\theta_0^{hk}\|_E^2\\
&
+ C k\sum_{j=1}^n \left(\|\dot{v}_j -\delta v_j\|_H^2 + \|v_j-w_j^h\|^2_V + \|u_j-u_j^{hk}\|^2_V+\|\dot{\theta}_j-\delta\theta_j\|_{L^2(\Omega)}^2 \right) \\
&
+Ck\sum_{j=1}^n \left( \|\theta_j-\eta_j^h\|_E^2+\|v_j-w_j^h\|_{L^2(\Gamma_C;\real^d)} + \|v_j-v_j^{hk}\|_{H}^2+ \|\theta_j-\theta_j^{hk}\|_{L^2(\Omega)}^2\right) \\[2mm] 
& 
+2k\sum_{j=1}^n ( \delta v_j -\delta v_j^{hk},v_j-w_j^h )_H
+2k \sum_{j=1}^n ( \delta\theta_j -\delta\theta_j^{hk},\theta_j-\eta_j^h )_{L^2(\Omega)}
\end{align*}
for all $\{w_j^h\}_{j=1}^n\subset V^h$ and $\{\eta_j^h\}_{j=1}^n\subset E^h$. 
\medskip
We have
\begin{align*}
&
\sum_{j=1}^n k ( \delta v_j-\delta v_j^{hk},v_j-w_j^h )_H \le \varepsilon \|v_n-v_n^{hk}\|_H^2 + C \|v_n-w_n^h\|_H^2 + C \|v_0-v_0^{hk}\|_H^2 \\
&
+ C\|v_1-w_1^h\|_H^2 + \sum_{j=1}^{n-1} 4k \|v_j-v_j^{hk}\|_H^2 + \frac{1}{k}\sum_{j=1}^{n-1} \|(v_j-w_j^h) - (v_{j+1} - w_{j+1}^h)\|_H^2
\end{align*}
and
\begin{align*}
&
\sum_{j=1}^n k ( \delta \theta_j-\delta \theta_j^{hk},\theta_j-\eta_j^h )_{L^2(\Omega)} \le \varepsilon \|\theta_n-\theta_n^{hk}\|_{L^2(\Omega)}^2 + C \|\theta_n-\eta_n^h\|_{L^2(\Omega)}^2  \\
&
+ C \|\theta_0-\theta_0^{hk}\|_{L^2(\Omega)}^2 + C\|\theta_1-\eta_1^h\|_{L^2(\Omega)}^2 + \sum_{j=1}^{n-1} 4k \|\theta_j-\theta_j^{hk}\|_{L^2(\Omega)}^2 \\
& + \frac{1}{k}\sum_{j=1}^{n-1} \|(\theta_j-\eta_j^h) - (\theta_{j+1} - \eta_{j+1}^h)\|_{L^2(\Omega)}^2.
\end{align*}
Next, we estimate
\begin{equation}
\|u_j-u_j^{hk}\|_V \le \|u_0-u_0^{hk}\|_V + \sum_{l=1}^j k\|v_l-v_l^{hk}\|_V + I_j,\label{KB_1}
\end{equation}
where $I_j$ is the integration error
\begin{equation*}
I_j = \left\|\int_0^{t_j} v(s)\,ds -\sum_{l=1}^j k v_l\right\|_V.
\end{equation*}
Since $I_j\le k\|u\|_{H^2(0,T;V)}$, we get
\begin{equation}
\|u_j-u_j^{hk}\|_V^2 \le C\left(\|u_0-u_0^{hk}\|^2_V + j\sum_{l=1}^j k^2 \|v_l-v_l^{hk}\|_V^2 + k^2\|u\|^2_{H^2(0,T;V)}\right)     \nonumber
\end{equation}
Using the fact that $Nk=T$, we have the estimate
\begin{align}\label{KB1}
&
\sum_{j=1}^n k\|u_j-u_j^{hk}\|_V^2\le CT \left(\|u_0-u_0^{h}\|^2_V + k^2 \|u\|^2_{H^2(0,T;V)}\right)\nonumber \\
& + T\sum_{j=1}^n k\sum_{l=1}^j k \|v_l-v_l^{hk}\|_V^2.    
\end{align}
\noindent We denote 
\[e_n=\|v_n-v_n^{hk}\|_H^2 + \sum_{j=1}^n k\|v_j-v_j^{hk}\|_V^2 + \|\theta_n - \theta_n^{hk}\|_{L^2(\Omega)}^2 + \sum_{j=1}^n k\|\theta_j-\theta_j^{hk}\|_E^2,\]
and
\begin{align*}
&
g_n =  k\sum_{j=1}^n \left(\|\dot{v}_j-\delta v_j\|_H^2 + \|v_j-w_j^h\|^2_V+\|\dot{\theta}_j -\delta \theta_j\|_{L^2(\Omega)}^2 + \|\theta_j -\eta_j^h\|^2_E\right)\\
&
+ \frac{1}{k}\sum_{j=1}^{n-1} \left(\|v_j-w_j^h - (v_{j+1} -w_{j+1}^h)\|_H^2+\|\theta_j -\eta_j^h - (\theta_{j+1}-\eta_{j+1}^{h})\|^2_{L^2(\Omega)}\right) \\
&
+k \sum_{j=1}^n \|v_j -w_j^{h}\|_{L^2(\Gamma_C;\real^d)}  
+ k^2\|u\|^2_{H^2(0,T;V)}+\|v_0-v_0^{hk}\|_V^2   \\
& +\|v_0-v_0^{hk}\|_H^2+ \|v_1-u_0^h\|_V^2+ \|\theta_0-\theta_0^{hk}\|_{L^2(\Omega)}^2+ \|\theta_0 -\theta_0^{hk}\|_E^2 \\[2mm]
&+\|\theta_1-\eta_1^h\|^2_{L^2(\Omega)}  + \|v_1 -w_1^h\|_H^2+ \|v_n-w_n^h\|_H^2+\|\theta_n-\eta_n^h\|^2_{L^2(\Omega)}.
\end{align*}
Then, we find that
\begin{equation}
e_n\le Cg_n + C\sum_{j=1}^n k e_j \ \ \mbox{for} \ n=1,\ldots,N
\end{equation}
with $C>0$. 
Using Lemma~\ref{gronwall}, we obtain the desired estimate (\ref{KB2}).
$\hfill{\Box}$

Now we consider a special situation in which the spaces $E^h$ and $V^h$ are based on the piecewise affine finite element approximation. Assume that $\Omega$ is a polygonal domain in $\real^2$  or a polyhedral domain in $\real^3$. Then each of parts of the boundary
$\Gamma_N$ and $\Gamma_C$ is a sum of segments (in the case $\Omega\subset\real^2$) or polygons (if $\Omega\subset\real^3$),
having mutually disjoint interiors i.e.
\[ \bar{\Gamma}_j=\bigcup_{i=1}^{N_j}\Gamma_{j,i},\quad j\in\{N,C\}.\]
Let  $\{\T^h\}$ be a family of regular triangulations of the set $\bar{\Omega}$ of diameter $h$ (cf. \cite{Ciarlet}, Chapter 3,  {\S} 3.1) which introduces a division of $\bar{\Omega}$ into triangles/tetrahedrons compatible with the division of the boundary $\partial\Omega$ into parts $\Gamma_{j,i}$, $1\le i\le N_j$, $j\in\{N,C\}$, in that sense that if any edge of a triangle/any face of tetrahedron that belongs to $\T^h$ has an intersection of positive boundary measure with any of sets $\Gamma_{j,i}$, then the whole edge/face is contained in  $\Gamma_{j,i}$. For an element $T\in \T^h$, let
$P_1(T)$ denote the space of first-order polynomials on $T$. Then $E^h$ and $V^h$ are defined as the spaces of continuous piecewise affine functions on elements $T\in T^h$, that is, 
\begin{align}
&E^h=\{w^h\in C(\bar{\Omega})\,\,\,|\,\, \,w^h|_T\in [P_1(T)]
\,\,\,\,\,{\rm for\,\, all}\,\,T\in {\T}^h\}\subset E.
\label{Eh}\\
&V^h=\{w^h\in C(\bar{\Omega};\real^d)\,\,\,|\,\, \,w^h|_T\in [P_1(T)]^d
\,\,\,\,\,{\rm for\,\, all}\,\,T\in {\T}^h\}\subset V.
\label{Vh}
\end{align}

\begin{Corollary} \label{col:conv} Suppose that the solution $(u,\theta)$ of Problem ${\cal P_E}$ satisfies the following regularity assumptions
	\begin{align*}
	&
	u\in C^1(0,T;H^2(\Omega;\real^d))\cap H^3(0,T;H)\cap H^2(0,T;V),\,\, \dot{u}_\tau\in C(0,T;H^2(\Gamma_C;\real^d)), \  \\
	&
  \theta\in C(0,T;H^2(\Omega))\cap H^2(0,T;L^2(\Omega))\cap H^1(0,T;E).
	\end{align*}
	Moreover, we assume that the initial data satisfy
	\begin{align*}
	u_0\in H^2(\Omega;\real^d),\,\, v_0\in H^1(\Omega;\real^d),\,\, \theta_0\in H^2(\Omega).
	\end{align*}
	Then, for $h<1$, we have the optimal order error estimate
    \begin{align}\label{error2}
        \max_{1\le n\le N} \{\|u_n-u_n^{hk}\|_V^2 + \|v_n-v_n^{hk}\|_H^2 + \|\theta_n-\theta_n^{hk}\|_{L^2(\Omega)}^2\}
        \nonumber \\
        + k\sum_{j=1}^N \|v_j-v_j^{hk}\|_V^2+ k\sum_{j=1}^N \|\theta_j-\theta_j^{hk}\|_E^2 \le C(h+k).
    \end{align}
\end{Corollary}
{\bf Proof.} To prove (\ref{error2}), we will use Theorem \ref{main2} and estimate each of the terms on the right-hand side of (\ref{KB2}). We start with the initial conditions. Let us define elements $u_0^h=\Pi^hu_0$, $v_0^h=\Pi^hv_0$, $\theta_0^h=\Pi^h\theta_0$ as finite element interpolants of $u_0$, $v_0$ and $\theta_0$, respectively (see \cite{Ciarlet}, (2.3.29)) By the standard finite element interpolation error estimates (\cite{AH2009, BS, Ciarlet}) we have the following
\begin{align*}
&\|u_0-u^h_0\|_V\le ch\,\|u_0\|_{H^2(\Omega;\mathbb{R}^d)},\\
&\|v_0-u^h_1\|_H\le ch\,\|v_0\|_{H^1(\Omega;\mathbb{R}^d)},\\
&\|\theta_0-\theta_0^h\|_E \le ch \|\theta_0\|_{H^2(\Omega)}.
\end{align*}
Let $w^h(t)=\Pi^h\dot{u}(t)\in V^h$ be the finite element interpolant of $\dot{u}(t)$, 
for $t\in (0,T]$.  Note that $w_\tau^h(t)=(\Pi^h\dot{u}(t))_\tau$ is the continuous piecewise-linear interpolant of $\dot{u}_\tau(t)$ on $\Gamma_C$.  Similarly, let $\eta^h(t)=\Pi^h\theta(t)\in E^h$ be the finite element interpolant of $\theta(t)$, a.e. $t\in [0,T]$. Then we have
\begin{align*}
&\|\dot{u}(t)-w^h(t)\|_V\leq ch\|\dot{u}(t)\|_{H^2(\Omega;\real^d)},\\
&\|\dot{u}(t)-w^h(t)\|_H \le ch^2 \|\dot{u}(t)\|_{H^2(\Omega;\real^d)},\\
& \|\dot{u}_\tau(t)-w   _\tau^h(t)\|_{L^2(\Gamma_C;\real^d)} \leq ch^2\|\dot{u}_\tau(t)\|_{H^2(\Gamma_C;\real^d)},\\
&\|\theta(t)-\eta^h(t)\|_E \le ch \|\theta(t)\|_{H^2(\Omega)},\\
&\|\theta(t)-\eta^h(t)\|_{L^2(\Omega)} \le ch^2 \|\theta(t)\|_{H^2(\Omega)}.
\end{align*}

Applying (\ref{KB2}) with $w_n^h:=w^h(t_n)$ and $\eta_n^h:=\eta^h(t_n)$, $n=1,\dots,N$, we get the following estimates
\begin{align*}
&k\sum_{j=1}^N \|v_j-w_j^h\|^2_V \le ch^2\|\dot{u}\|^2_{C(0,T;H^2(\Omega;\real^d))},
\\[2mm]
&\max_{1\le n\le N} \|v_{n\tau} - w_{n\tau}^h\|_{L^2(\Gamma_C;\real^d)}\leq ch^2\|\dot{u}_\tau\|_{C(0,T;{H^2(\Gamma_C;\real^d))}},\\[2mm]
&\max_{1\le n\le N} \|v_n-w_n^h\|_H^2\leq ch^4\|\dot{u}\|^2_{C(0,T;H^2(\Omega;\real^d))},\\[2mm]
& k\sum_{j=1}^N \|\theta_j - \eta_j^h\|_E^2\le ch^2\|\theta\|^2_{C(0,T;H^2(\Omega))},\\[2mm]
&\max_{1\le n\le N} \|\theta_n-\eta_n^h\|^2_{L^2(\Omega)}\leq ch^4\|\theta\|^2_{C(0,T;H^2(\Omega))}.
\end{align*}
Moreover, similarly to Section 5 of \cite{HSS}, we get
\begin{equation*}
k\sum_{j=1}^N \|\dot{v}_j -\delta v_j\|_H^2 \le ck^2 \|u\|_{H^3(0,T;H)}^2,
\end{equation*}
\begin{equation*}
k\sum_{j=1}^N \|\dot{\theta}_j -\delta \theta_j\|_{L^2(\Omega)}^2 \le ck^2 \|\theta\|_{H^2(0,T;L^2(\Omega))}^2,
\end{equation*}
\begin{equation*}
\frac{1}{k}\sum_{j=1}^{N-1} \|(v_j -w_j^h) - (v_{j+1}-w_{j+1}^h)\|_H^2 \le ch^2\|u\|_{H^2(0,T;V)}^2,
\end{equation*}
\begin{equation*}
\frac{1}{k}\sum_{j=1}^{N-1} \|(\theta_j -\eta_j^h) - (\theta_{j+1}-\eta_{j+1}^h)\|_{L^2(\Omega)}^2 \le ch^2\|\theta\|_{H^1(0,T;E)}^2.
\end{equation*}
Then the error bound (\ref{error2}) follows from (\ref{KB2}) combined with (\ref{KB_1}). The proof of the corollary is complete.
$\hfill{\Box}$


\section{Numerical simulations}\label{SecNum}
In this section, we present the results of our numerical simulations. Such simulations are conducted to compare the empirical and theoretical convergence results of the proposed fully discrete scheme and visualize the role of various factors in the considered model. In Section \ref{sec:sim:2d}, an academic two-dimensional sample problem is solved to obtain empirical error estimates and compare them with the theoretical estimates presented in Corollary \ref{col:conv}. Next, in Section \ref{sec:sim:3d} we consider the three-dimensional problem with a more complex shape of the body.

We employ FEM and use the definition of spaces $E^h$ and $W^h$ presented in (\ref{Eh}) and (\ref{Vh}), respectively. For the sake of clarity, we assume dimensionless models with mass density and heat capacity $\rho = c_p = 1$. In our simulations, we use the original open source package \textit{conmech} (\url{https://github.com/KOS-UJ/conmech}). The code is written in Python and uses just-in-time compilation mechanisms provided by packages \textit{numba} \cite{NUMBA} and \textit{JAX} \cite{JAX}.

For two-dimensional simulations, we use the direct optimization method with the Schur complement method to decrease the dimension of the system \cite{OJB}.
We use Powell's conjugate direction method from the \textit{SciPy} \cite{SCIPY} package 
since the minimized functional is not necessarily differentiable. 

In the case of three-dimensional simulations, to speed up the calculations, we perform computations on GPU using the L-BFGS solver implemented using \textit{JAX}. We use \textit{TetWild} \cite{TETWILD} to generate the volumetric mesh from the surface mesh obtained from the Stanford 3D Scanning Repository. The results are visualized using \textit{Blender} (\cite{BLENDER}).

\subsection{Two-dimensional numerical example} \label{sec:sim:2d}

We consider an isotropic homogeneous domain $\Omega = [0, 1.5] \times [0, 1]$, set $T=1$, and assume the following boundary conditions
\begin{equation*}
\Gamma_N = (\{0\} \times (0, 1]) \cup ([0, 1.5] \times \{1\}) \cup (\{1\} \times (0, 1]), \quad \Gamma_C = [0, 1.5] \times \{0\}.
\end{equation*}

\noindent
The domain $\Omega$ represents the cross-section of a three-dimensional viscoelastic body.
We use the Kelvin-Voigt type short-memory viscoelastic law. 
The viscosity operator $\mathcal{A}$, elasticity operator $\mathcal{B}$, heat expansion tensor $C$ and the thermal conductivity operator $K$ are defined by
\begin{align*}
 \mathcal{A}(t, {\tau}) &= 2\phi{\tau} + \xi tr({\tau})I,\quad {\tau} \in \mathbb{S}^3,\, t \in [0, T],\\ 
 \mathcal{B}({\tau}) &= 2\mu{\tau} + \lambda tr({\tau})I,\quad {\tau} \in \mathbb{S}^3,\\ 
 C &= c\, I,\\ 
 K(t, {\zeta}) &= \kappa\, {\zeta},\quad {\zeta} \in \mathbb{R}^3,\, t \in [0, T].
\end{align*}

\begin{figure}[ht!]
\centering
    \includegraphics[width=0.8\linewidth]{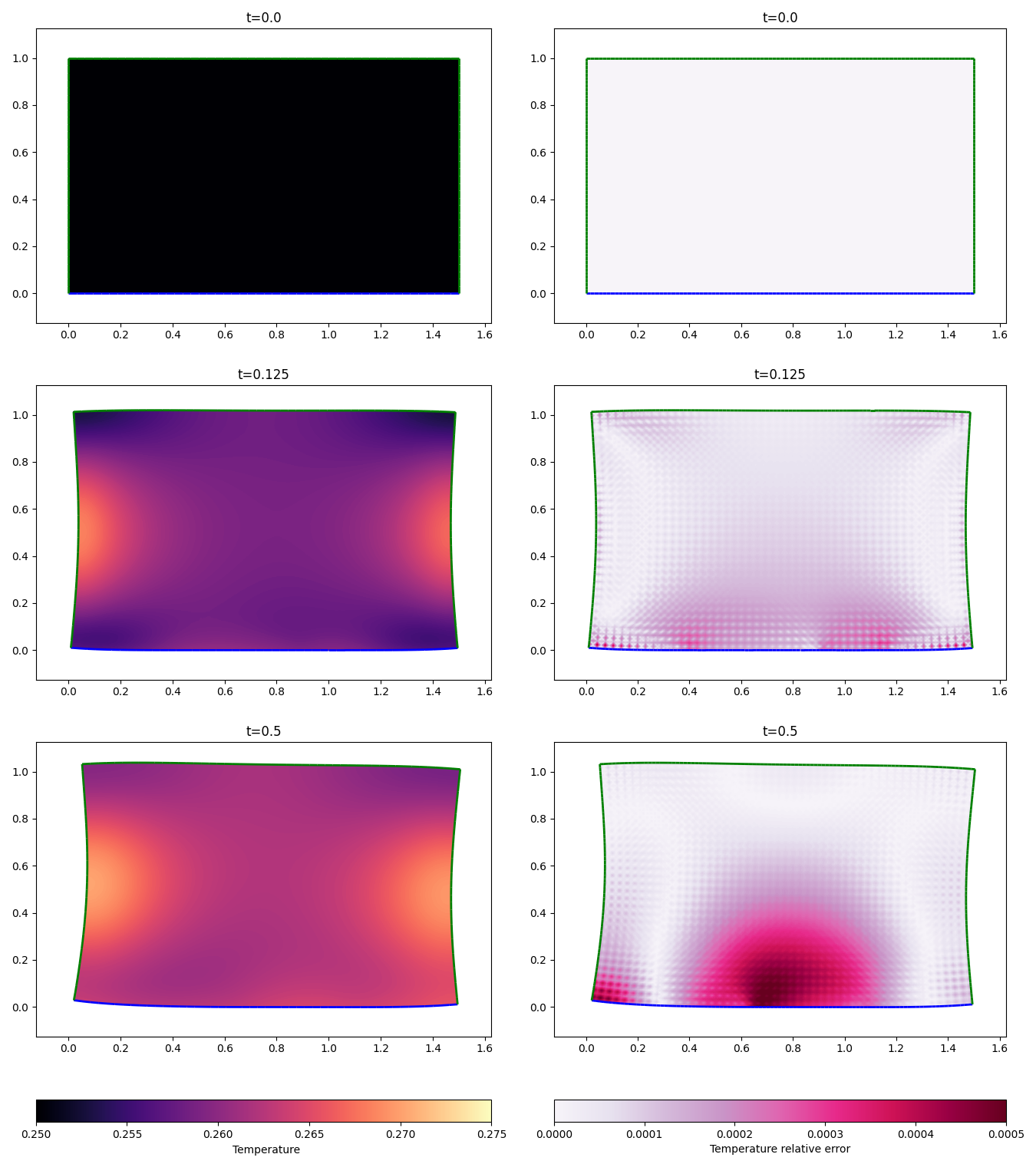}
    \caption{Solution at chosen time moments $t$.  Plots present the temperature on the left and relative temperature error on the right for $h=2^{-5}$ and $k=2^{-9}$.} \label{fig:2d}
\end{figure}

\noindent
Here, $I \in \mathbb{R}^{2\times 2}$ is the identity matrix,  $\text{tr}(\cdot)$ is the trace of a square matrix, $\lambda$ and $\mu$ are the Lame coefficients, while $\xi$ and $\phi$ represent the viscosity coefficients, $\lambda, \mu, \xi, \phi, c, k >0$.
In the initial state, the body rests on the foundation, that is, the initial displacement $u_0$ and the initial velocity $v_0$ are equal to zero at any point in $\Omega$, and the gap between the body and the foundation $g_0 = 0$ at any point on the contact boundary $\Gamma_C$.
The viscoelastic properties of the body are characterized by $\phi=4.5, \xi=10.5, \nu=45, \lambda = 105$.
The thermal effect is parameterized by the heat expansion coefficient $c = 0.5$ and the thermal conductivity coefficient $\kappa = 0.033$.
The internal heat source $g$ and the heat exchange between the body and the foundation $r$ are modeled by the functions
\begin{align*}
g(x,t) = 0, \quad & x \in \Omega,\ t \in [0,T], \\
r(x, \theta) = r_F(\theta - \theta_F), \quad & x \in \Gamma_C,\ \theta \in \mathbb{R},
\end{align*}
where $r_F = 0.1$ and $\theta_F = 0.27$ denote the heat exchange coefficient and the constant foundation temperature, respectively.
We assume the initial temperature $\theta_0 = 0.25$, constant density of the volume forces $f_0 = (0, \,-1)^T$ inside the body, and non-zero traction force applying squeezing from both sides of the body, namely
\begin{equation*}
f_N(x_1, x_2) =
    \begin{cases}
        (48(0.25 - (x_2 - 0.5)^2), \, 0)^T, & \text{for} \,\, (x_1, x_2) \in (\{0\} \times (0, 1]),\\
        (-44(0.25 - (x_2 - 0.5)^2), \, 0)^T, & \text{for} \,\, (x_1, x_2) \in (\{1\} \times (0, 1]),\\
        (0, \,0)^T, &  \text{for} \,\, (x_1, x_2) \in ([0, 1.5] \times \{1\}).
    \end{cases} 
\end{equation*}

\noindent
On the contact boundary $\Gamma_C$ we set
\begin{align*}
    p_\nu({x}, \zeta) &= \left \{ 
    \begin{array}{lll}
        0, & \text{for} \,\, \zeta \in (-\infty,\, 0), & {x} \in \Gamma_C,\\
        10^3\, \zeta, & \text{for} \,\, \zeta \in [0,\, \infty), & {x} \in \Gamma_C,\\
    \end{array} \right.\\[2mm]
    j(\zeta) &= \log(\|\zeta\|+1) \quad \text{for} \,\, \zeta \in \mathbb{R}_{+}, \\[2mm]
    h(x, \eta) &= \left \{ 
    \begin{array}{lll}
        0, & \text{for} \,\, \eta \in (-\infty,\, 0), &{x} \in \Gamma_C, \\
        0.1 \eta, & \text{for} \,\, \eta \in [0,\, \infty), &{x} \in \Gamma_C.
    \end{array} \right.
\end{align*} 
as conditions governing the response in the normal direction, friction, and heat generated by friction. Note that $j(\cdot)$ is in the form of a non-differentiable and non-convex function.

Figure \ref{fig:2d} demonstrates the solution of the sample problem in which we squeeze the body that lies on the foundation.
The traction force from the left side is greater than from the right side, so we can observe slipping of the body on the foundation and therefore friction and heat generated by friction.
On the right side, we can see the corresponding relative temperature error computed in comparison with the reference solution with $h=2^{-7}$, $k=2^{-9}$ and with 148419 degrees of freedom.

\begin{figure}[ht!]
\centering
    \includegraphics[width=0.45\linewidth]{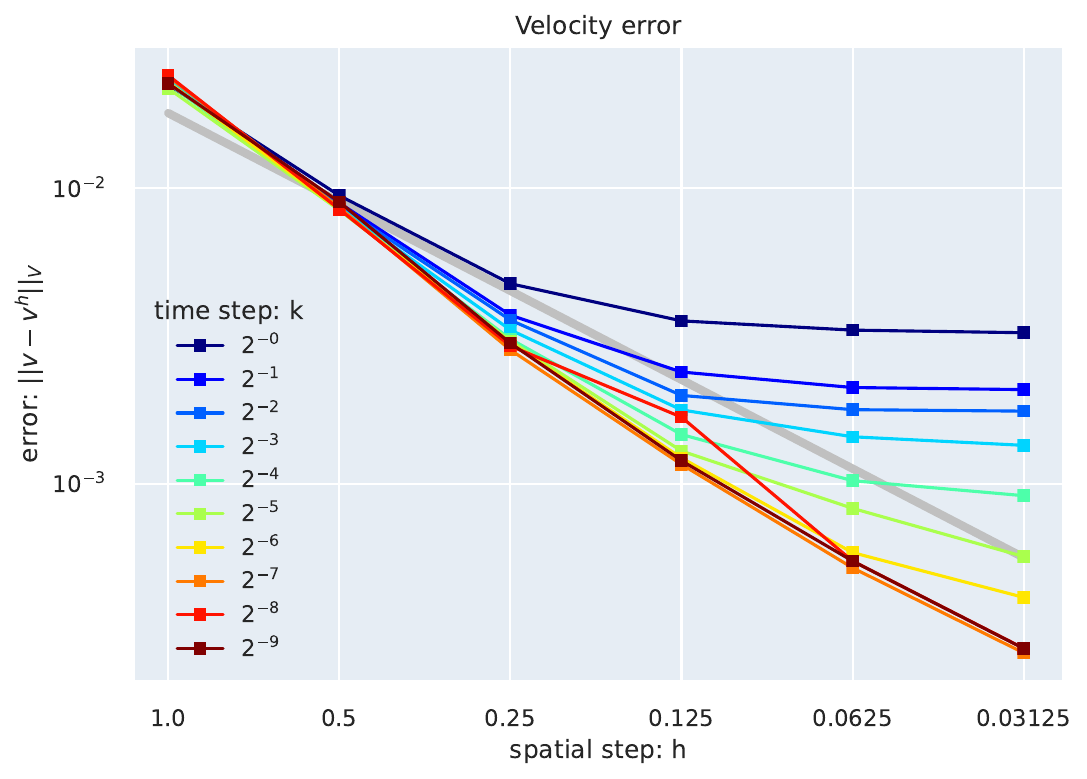} \quad
    \includegraphics[width=0.45\linewidth]{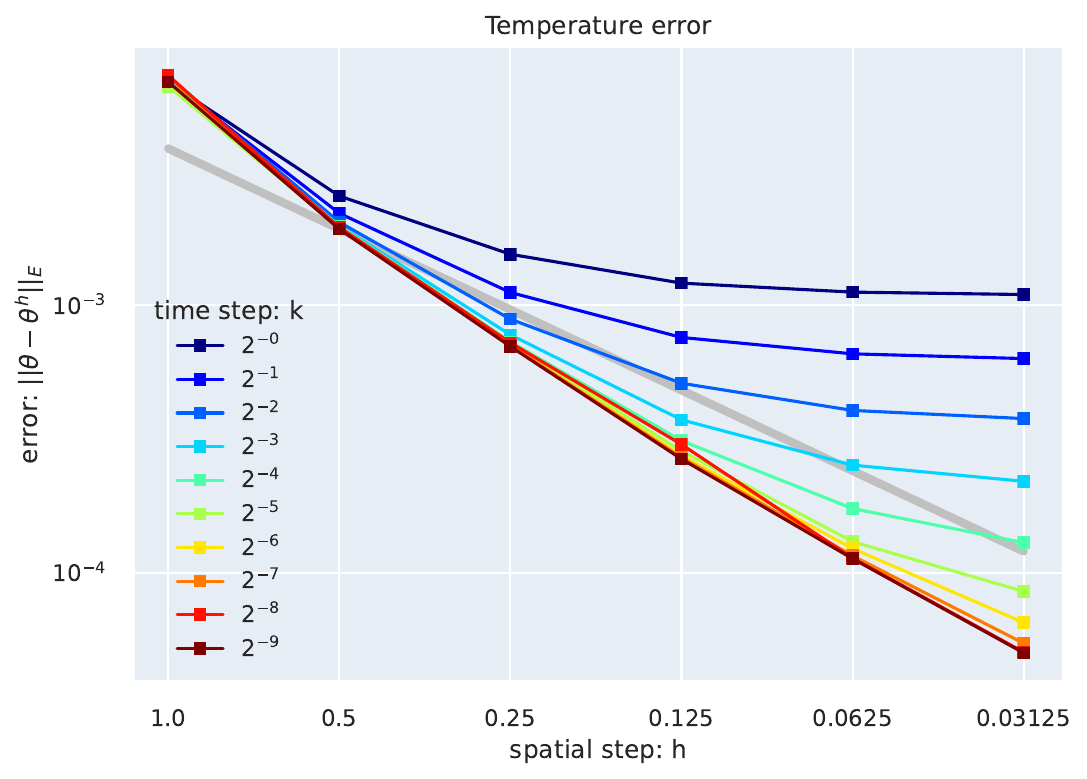}\\
    \caption{Numerical error estimate for velocity (left) and temperature (right).} \label{fig:err}
\end{figure}

Figure \ref{fig:err} depicts the empirical error estimation for the sample problem considered plotted on a log-log scale.
The color of the plot corresponds to the time step used to compute the solution. As we can see, for any fixed spatial step $h$ we can observe that the numerical error goes to zero as the time step goes to zero. The gray line shows the linear order of convergence for comparison with the theoretical results.

\subsection{Three-dimensional numerical example} \label{sec:sim:3d}

To visualize the influence of various modeled factors, we set $d=3$ and consider a volumetric mesh consisting of 42684 nodes and 230213 tetrahedrons. We set $\Gamma_{N} = \emptyset$ and $\Gamma_{C} = \Gamma$. The final time $T= 2.5$, the time step $k = 0.01$, and the results are presented at times $t \in \{0.5,\, 1,\, 1.5,\, 2,\, 2.5\}$. The position and shape of the foundation can be deduced from the attached illustrations. Body temperature is once again represented by its color, where blue is the coldest region and yellow is the hottest.
We apply analogous constitutive laws as in a two-dimensional example and set
\begin{align*}
 &\lambda = \mu = 10, \quad \xi = \phi = 20\\
 &c = 1, \quad \kappa = 0.1
\end{align*}
The initial displacement $u_0$ and the initial velocity $v_0$ are equal to zero at any point in $\Omega$. The internal force, the internal heat source and the heat exchange between the body and the foundation are modeled by
\begin{align*}
f_0(x,t) &= (-0.7, \,0, \,2)^T, \quad x \in \Omega,\ t \in [0,T], \\
g(x,t) &= 0, \quad x \in \Omega,\ t \in [0,T], \\
r(x, \theta) &= 0, \quad x \in \Gamma_C,\ \theta \in \mathbb{R}.
\end{align*}
On the contact boundary $\Gamma_C$ we take the following data
\begin{align*}
    p_\nu({x}, \zeta) &= \left \{ 
    \begin{array}{lll}
        0, & \text{for} \,\, \zeta \in (-\infty,\, 0), & {x} \in \Gamma_C,\\
        \hat{p}\, \zeta, & \text{for} \,\, \zeta \in [0,\, \infty), & {x} \in \Gamma_C,\\
    \end{array} \right.\\[2mm]
    j(\zeta) &= \hat{j}\, \log(\|\zeta\|+1) \quad \text{for} \,\, \zeta \in \mathbb{R}_{+}, \\[2mm]
    h(x, \eta) &= \left \{ 
    \begin{array}{lll}
        0, & \text{for} \,\, \eta \in (-\infty,\, 0), &{x} \in \Gamma_C, \\
        \hat{h}\, \eta & \text{for} \,\, \eta \in [0,\, \infty), &{x} \in \Gamma_C.
    \end{array} \right.
\end{align*} 
The values of the hardness of the foundation $\hat{p}$, the friction coefficient $\hat{j}$ and the coefficient of heat generated by friction $\hat{h}$ vary between examples.

\begin{figure}[ht!]
     \begin{center}
        \includegraphics[width=0.2\linewidth]{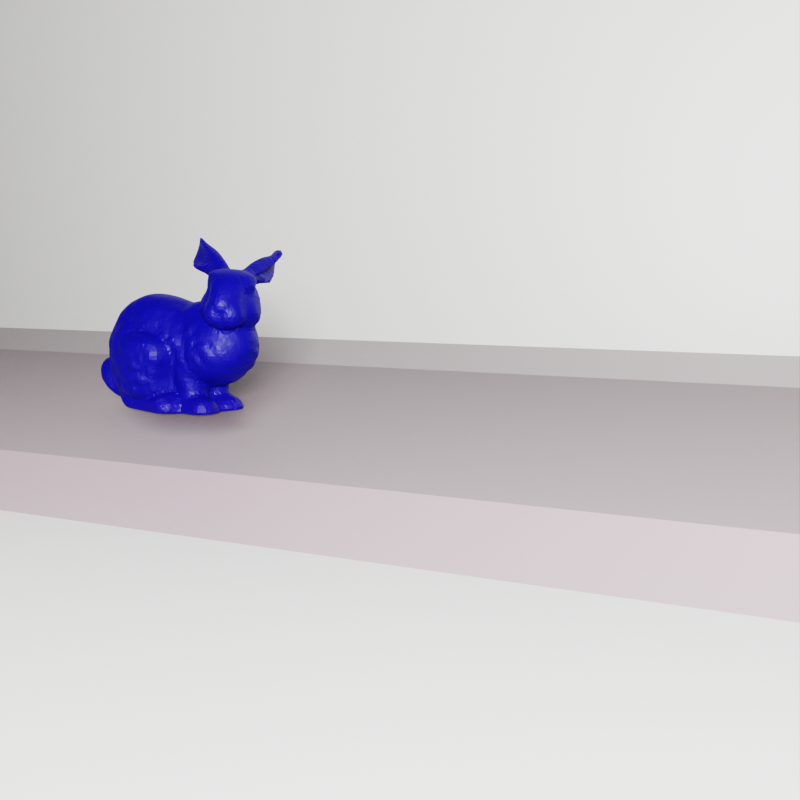}%
        \includegraphics[width=0.2\linewidth]{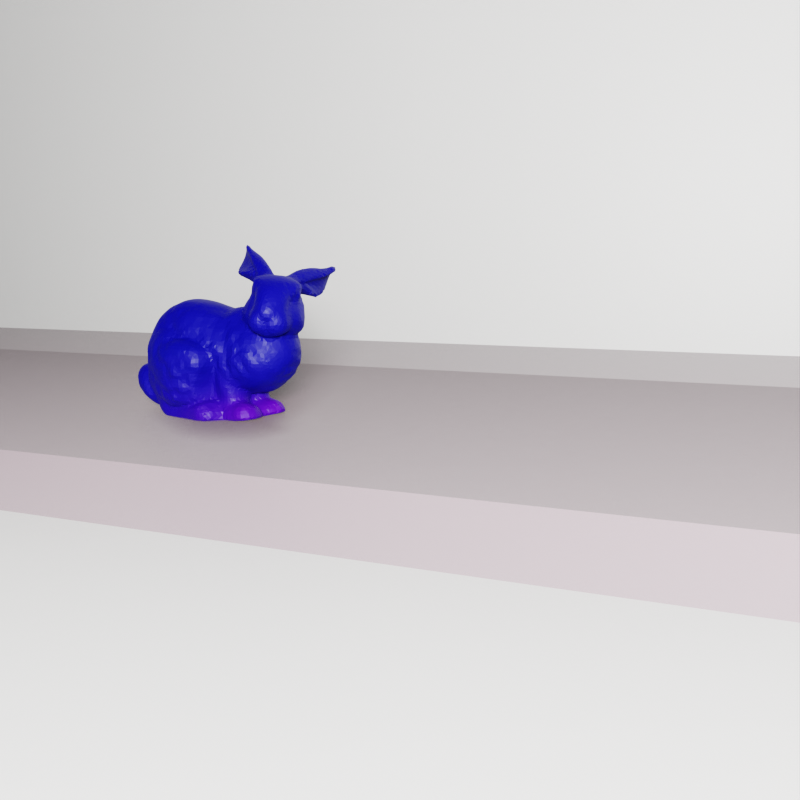}%
        \includegraphics[width=0.2\linewidth]{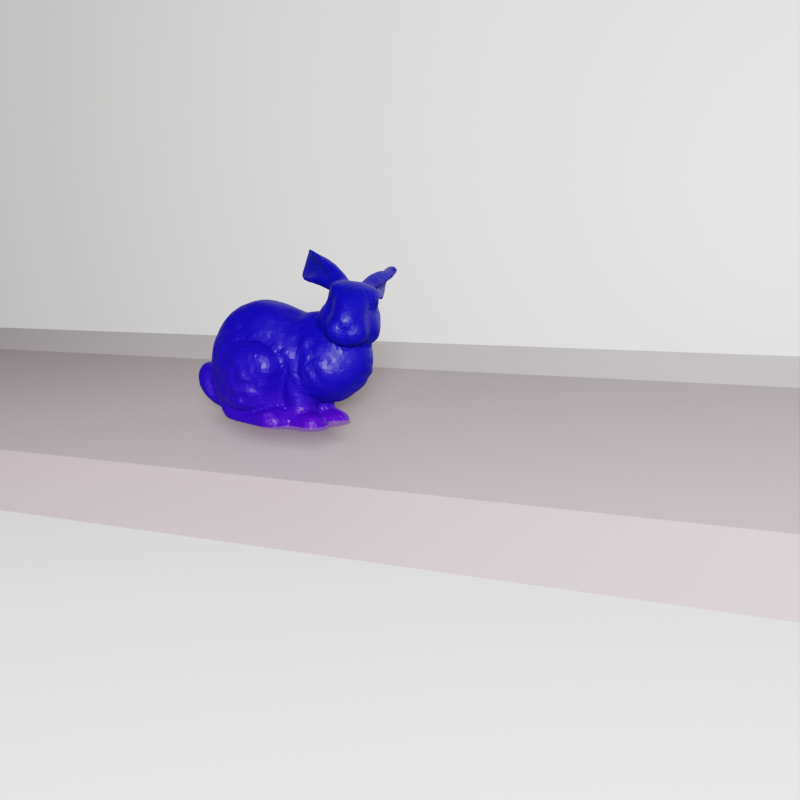}%
        \includegraphics[width=0.2\linewidth]{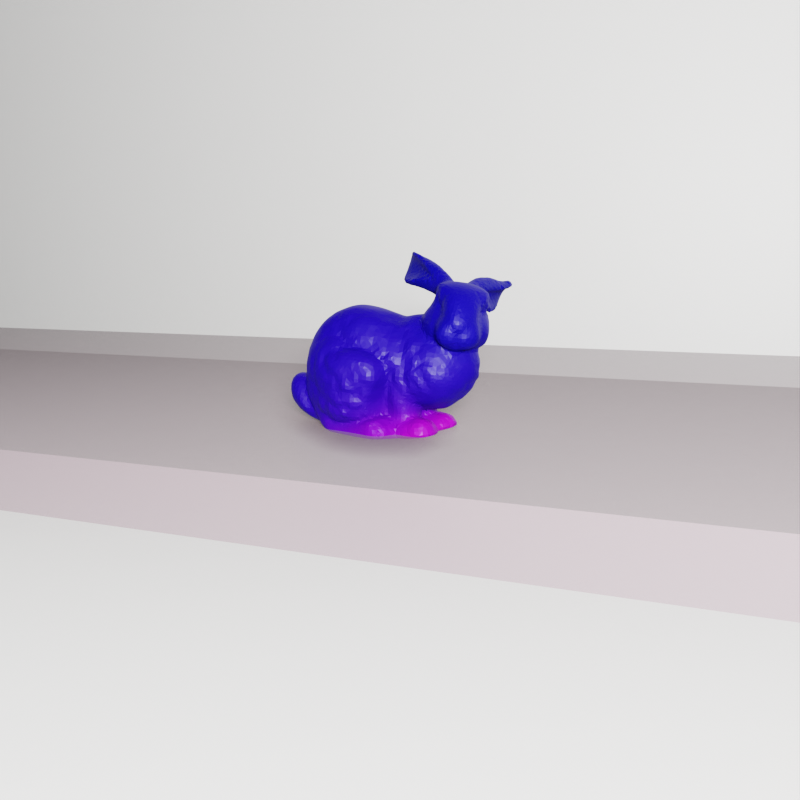}%
        \includegraphics[width=0.2\linewidth]{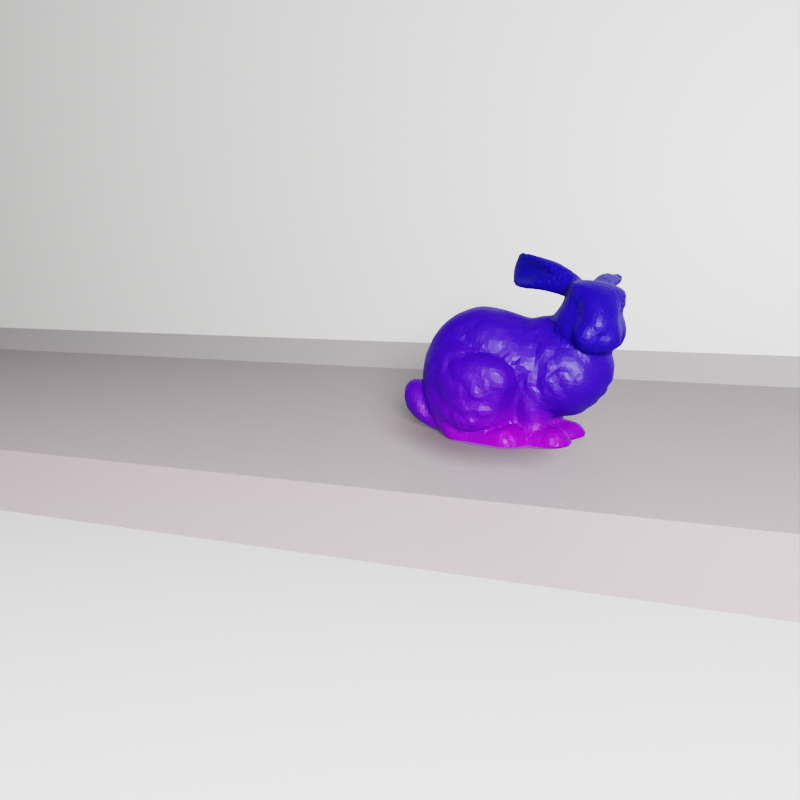}
        \caption{First simulation - base data
        } \label{fig:sim1}
        \hspace{5mm}
        \includegraphics[width=0.2\linewidth]{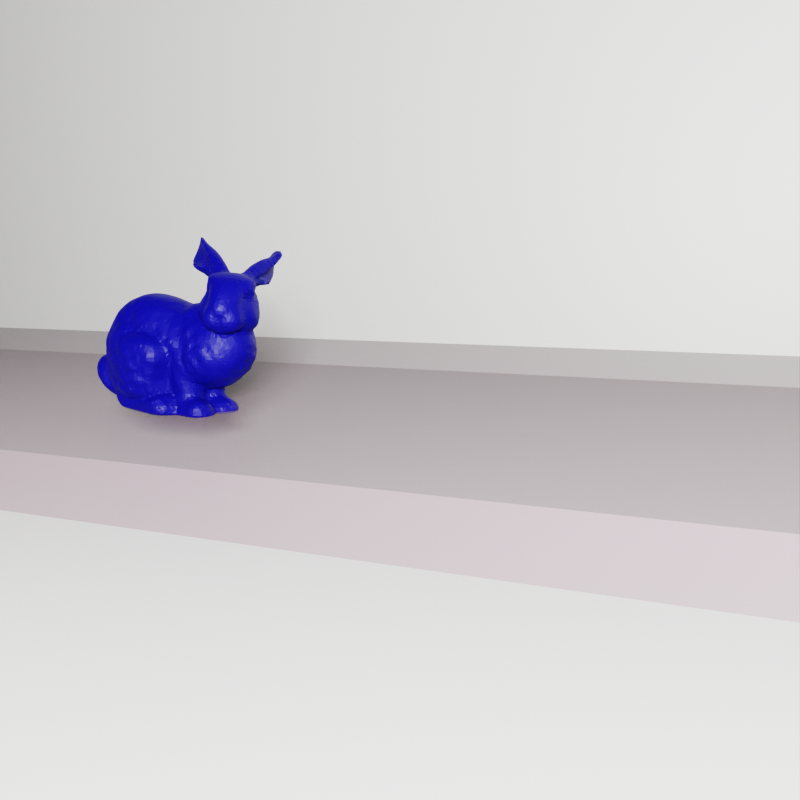}%
        \includegraphics[width=0.2\linewidth]{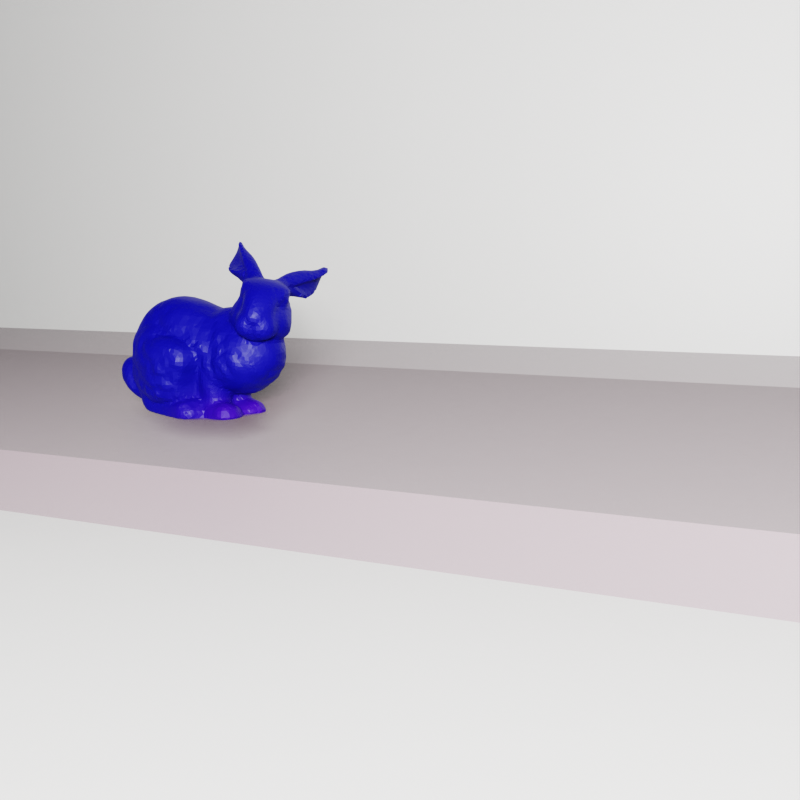}%
        \includegraphics[width=0.2\linewidth]{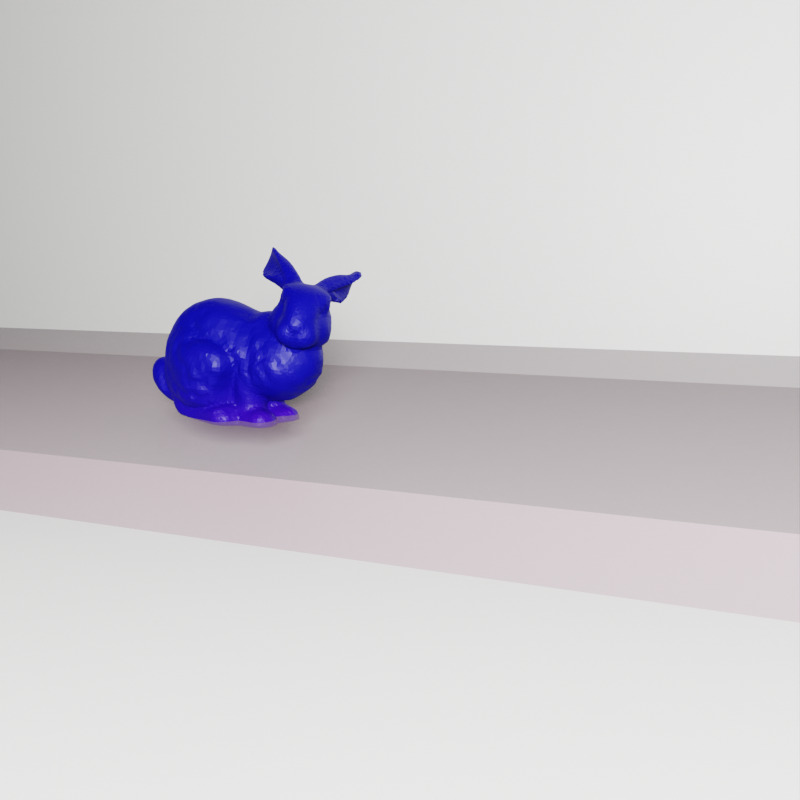}%
        \includegraphics[width=0.2\linewidth]{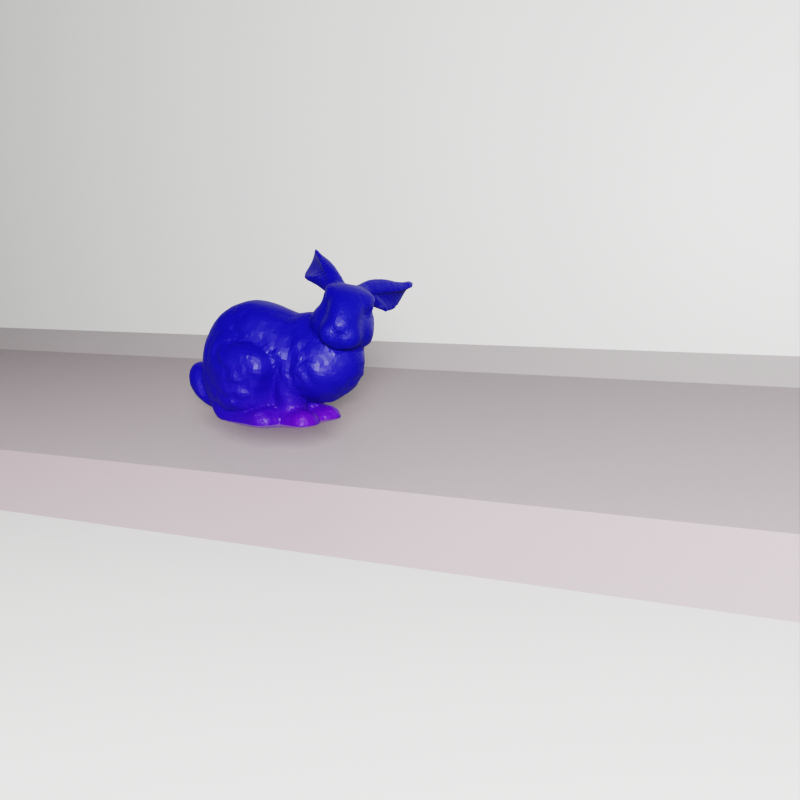}%
        \includegraphics[width=0.2\linewidth]{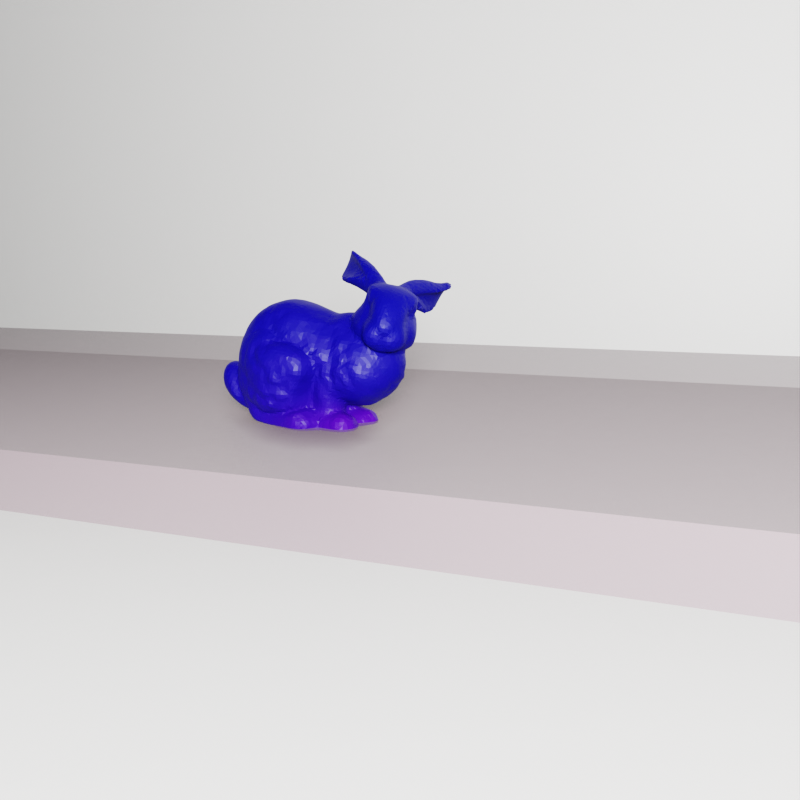}
        \caption{Second simulation - increased friction
        } \label{fig:sim2}
        \hspace{5mm}
        \includegraphics[width=0.2\linewidth]{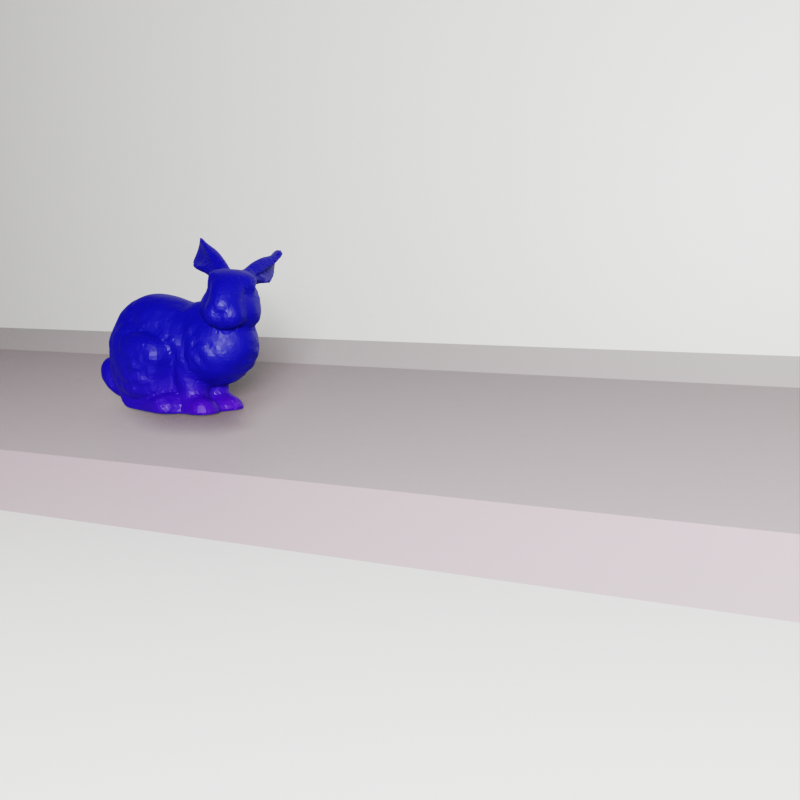}%
        \includegraphics[width=0.2\linewidth]{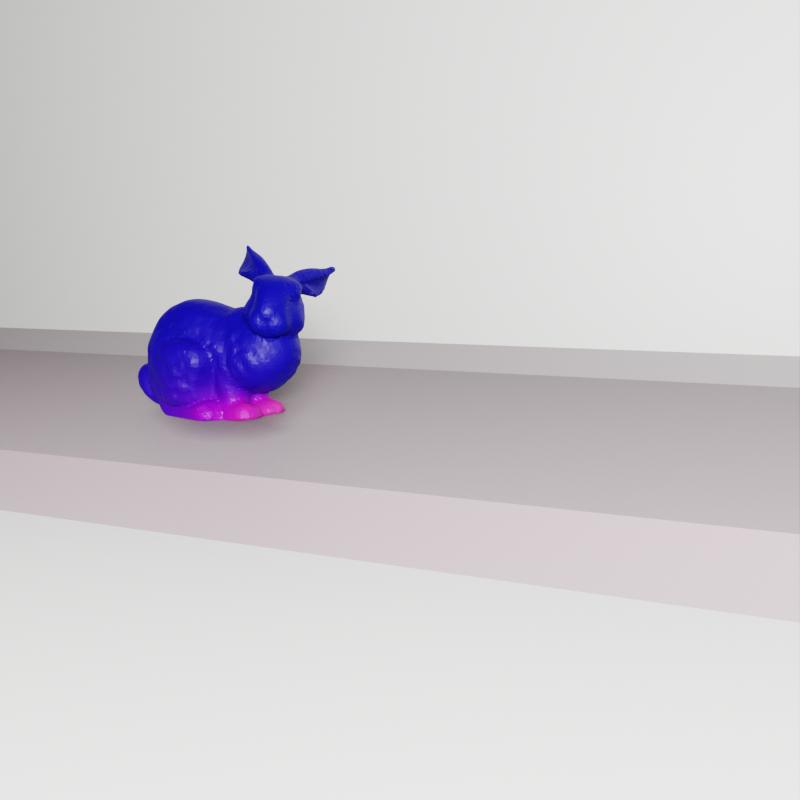}%
        \includegraphics[width=0.2\linewidth]{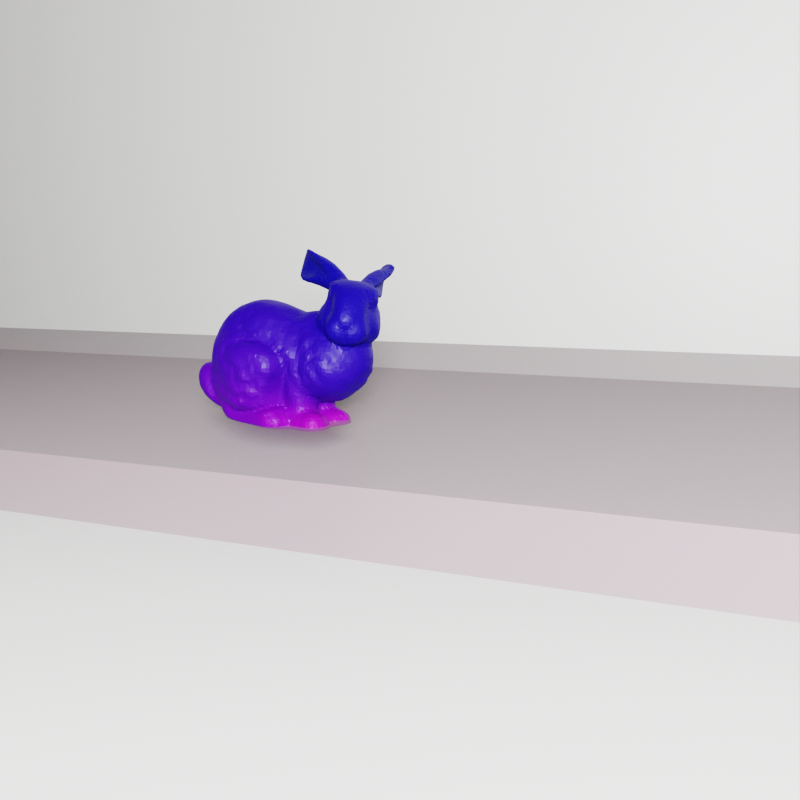}%
        \includegraphics[width=0.2\linewidth]{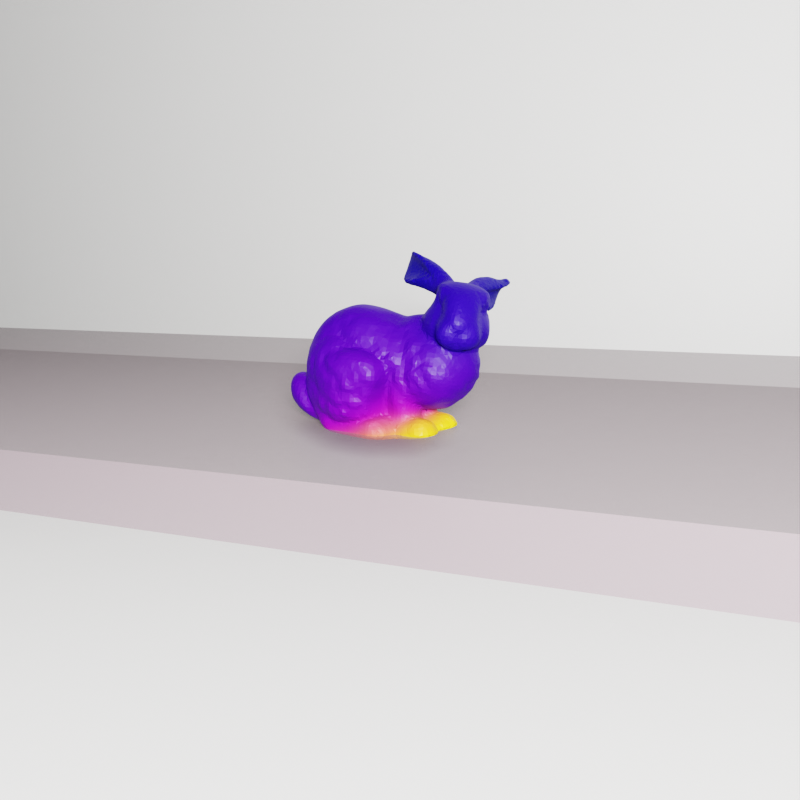}%
        \includegraphics[width=0.2\linewidth]{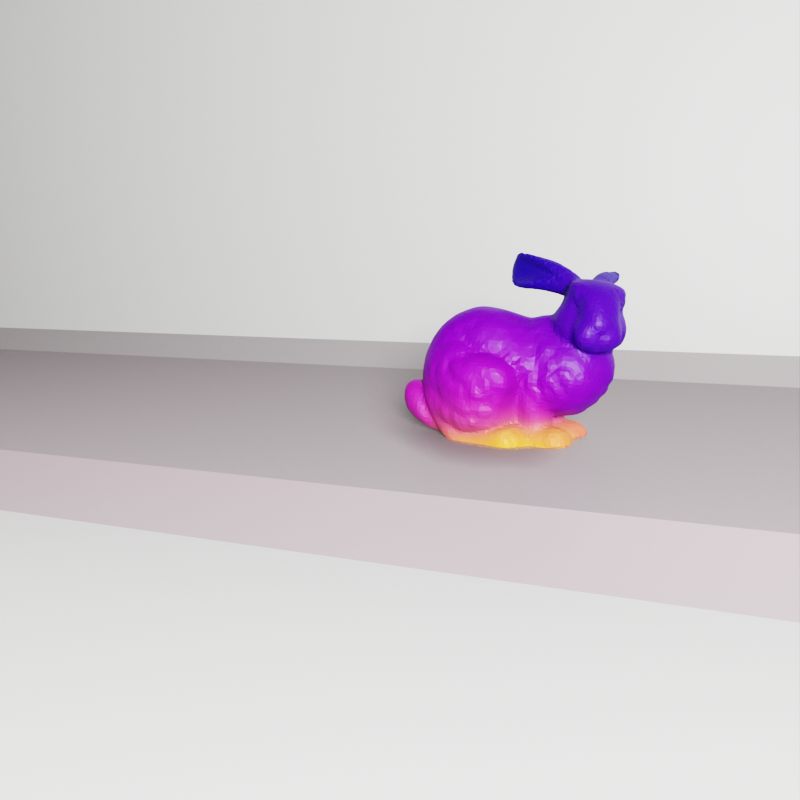}
        \caption{Third simulation - increased heat generated by friction
        } \label{fig:sim3}
        \hspace{5mm}
        \includegraphics[width=0.2\linewidth]{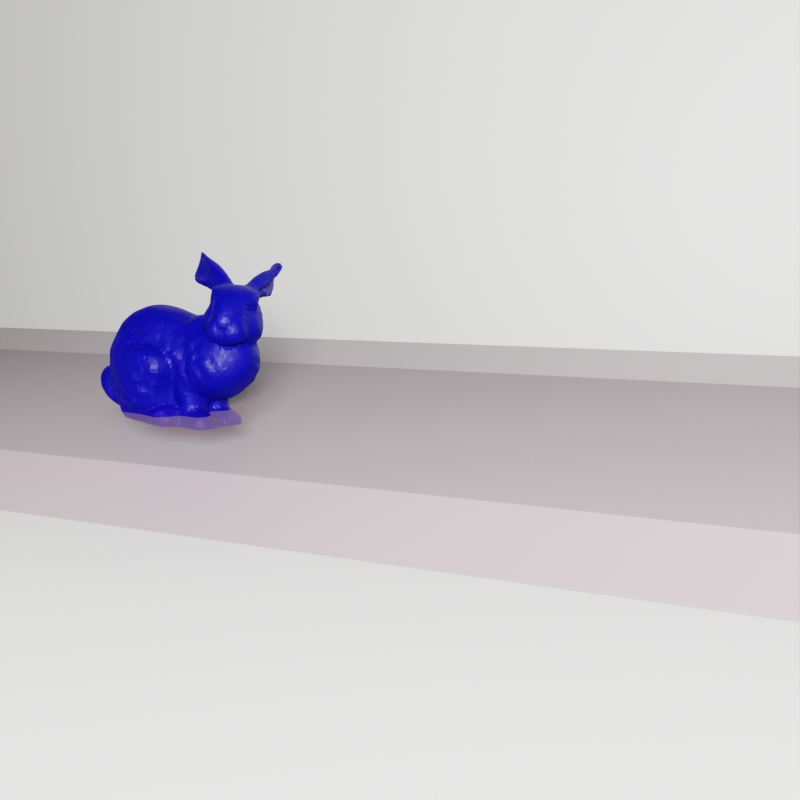}%
        \includegraphics[width=0.2\linewidth]{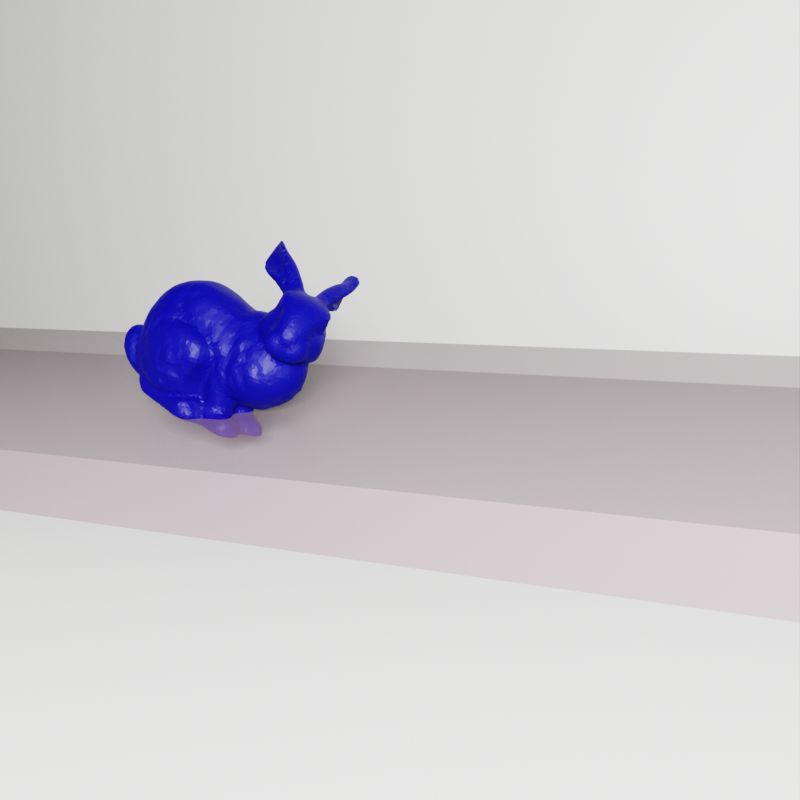}%
        \includegraphics[width=0.2\linewidth]{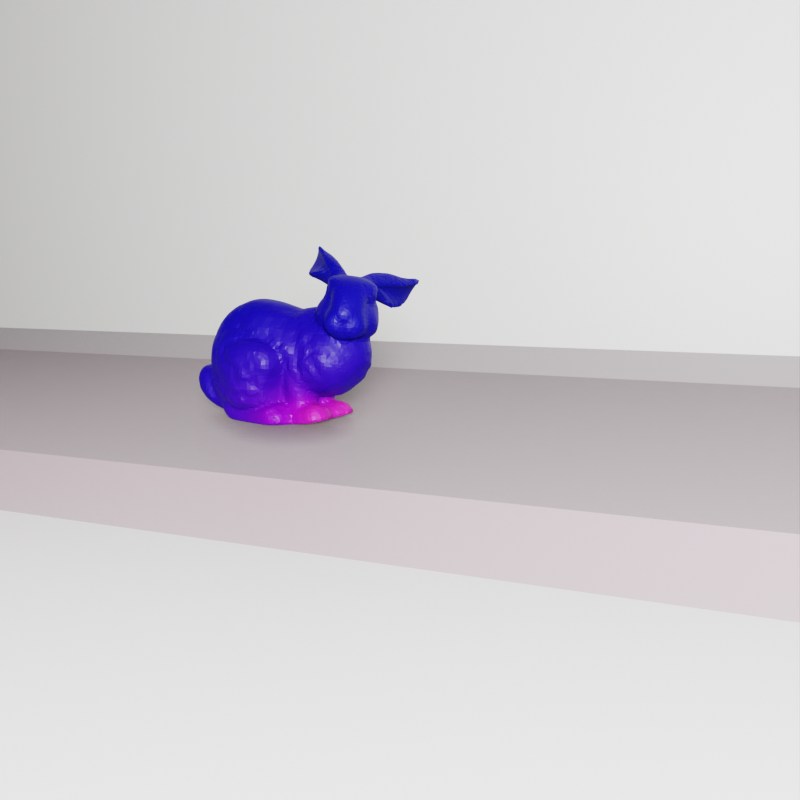}%
        \includegraphics[width=0.2\linewidth]{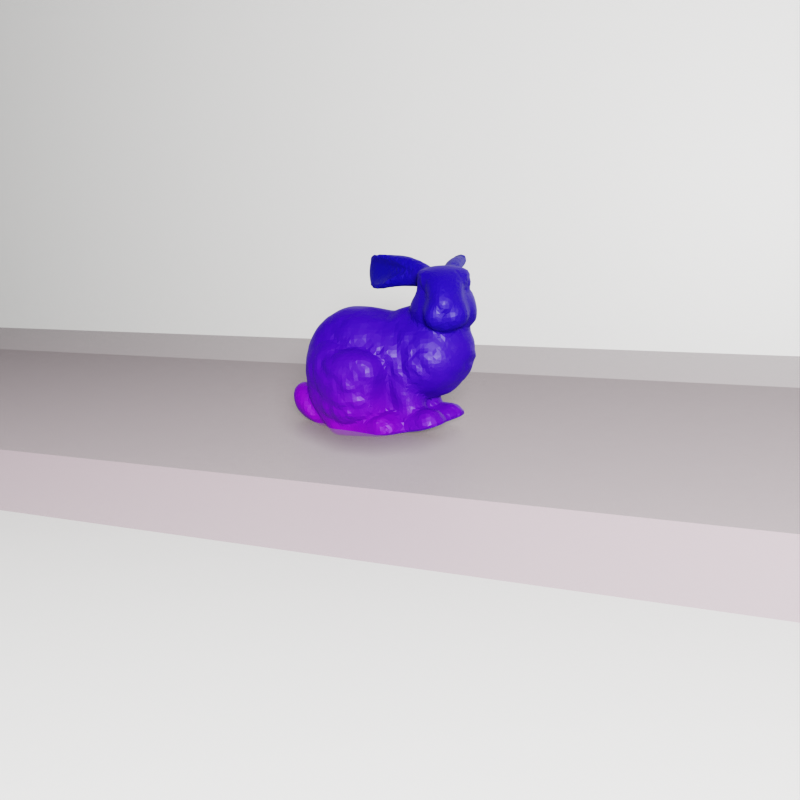}%
        \includegraphics[width=0.2\linewidth]{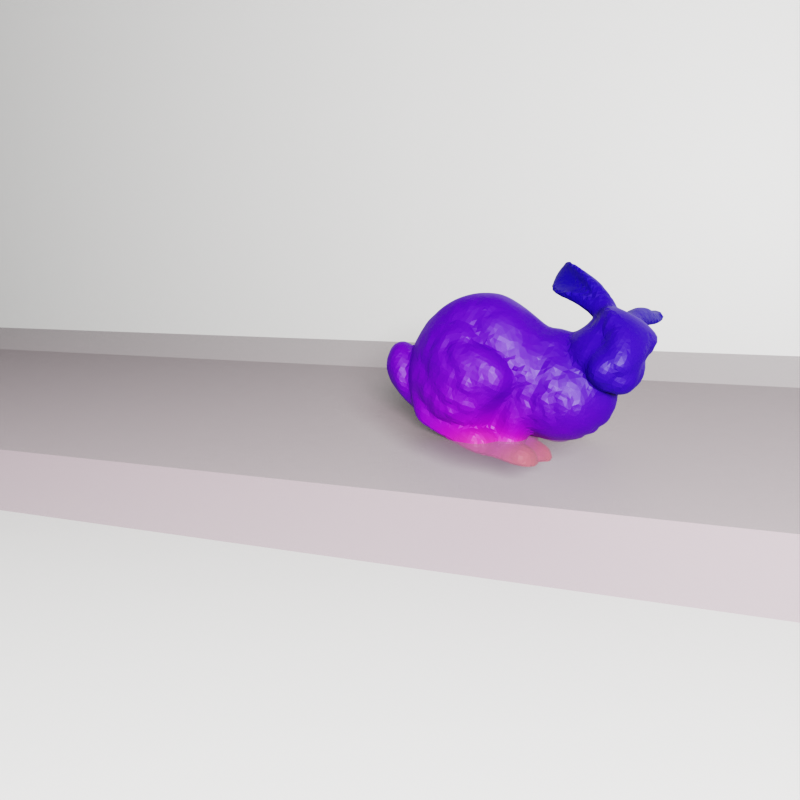}
        \caption{Fourth simulation - softer foundation
        } \label{fig:sim4}
      \end{center}     
\end{figure}

We now provide a brief overview of the results. The first simulation shows the data $\hat{p} = 200$, $\hat{j} = 0.1$ and $\hat{h} = 0.2$ and is presented in figure~\ref{fig:sim1}. This is the base configuration that will be modified in the following examples. We push the body down and to the right with force ${f}_0$. The body falls to the obstacle, rebounds, and slides to the right. We can also observe the dissipation of heat generated by friction.

The second simulation shows the data $\hat{p} = 200$, $\hat{j} = 1.0$ and $\hat{h} = 0.2$ and is presented in Figure~\ref{fig:sim2}. The only change compared to the base configuration is the increase in friction coefficient. This causes the body to cover a shorter distance and, as a result of the lower speed of the nodes at the contact boundary, the heat generated by friction is decreased.

The third simulation shows the data $\hat{p} = 200$, $\hat{j} = 0.1$ and $\hat{h} = 0.6$ and is presented in Figure~\ref{fig:sim3}. In this case, the coefficient of heat generated by friction is increased compared to the base configuration. This increases the heat generated at the contact boundary, which is then diffused inside the body.

The fourth simulation shows the data $\hat{p} = 50$, $\hat{j} = 0.1$ and $\hat{h} = 0.2$ and is presented in Figure~\ref{fig:sim4}. In this case, the boundary hardness coefficient is decreased. The softer foundation causes the body to initially dive further, and the resulting rebound produces a less stable trajectory and body rotation. We remark that the infinitesimal strain tensor considered in this paper results in a constitutive law that is not rotationally invariant. As a result, its application in the case of large body rotation would cause the volume of the body to increase unnaturally.

\bigskip
 
\section{Acknowledgements}\label{acknowledgement}
The research was
supported by the European Union's Horizon 2020 Research and Innovation Programme under
the Marie Sklodowska-Curie grant agreement No. 823731 CONMECH,
the Ministry of Science and Higher Education of Republic of Poland
under Grant Nos. 4004/GGPJII/H2020/2018/0 and 440328/PnH2/2019,
and the National Science Centre of Poland under Project No. 2021/41/B/ST1/01636.


\bibliographystyle{elsarticle-num}
\bibliography{sample}

\end{document}